%
\newif\ifloadreferences\loadreferencestrue
%
%
%
%
%
\let\myfrac=\frac%
\input eplain %
\let\frac=\myfrac%
\let\myfootnote=\footnote%
\input amstex \input epsf %
\let\footnote=\myfootnote%
%
%
\loadeufm\loadmsam\loadmsbm\message{symbol names}\UseAMSsymbols\message{,}%
\magnification 1200 %
\font\myfontdefault=cmr10%
\newif\ifmakebiblio%
\newif\ifinappendices%
\newif\ifundefinedreferences%
\newif\ifchangedreferences%
\makebibliofalse%
\undefinedreferencesfalse%
\changedreferencesfalse%
%
%
%
%
%
\def\setcatcodes{\catcode`\!=0 \catcode`\\=11}%
{\global\let\noe=\noexpand%
\catcode`\@=11 \catcode`\_=11 \setcatcodes%
!global!def!_@@internal@@makeref#1{%
!global!expandafter!def!csname #1ref!endcsname##1{%
!csname _@#1@##1!endcsname%
!expandafter!ifx!csname _@#1@##1!endcsname!relax%
    !write16{#1 ##1 not defined - run saving references}%
    !undefinedreferencestrue%
!fi}}%
!global!def!_@@internal@@makelabel#1{%
!global!expandafter!def!csname #1label!endcsname##1{%
!edef!temptoken{!csname #1info!endcsname}%
!ifloadreferences%
!expandafter!ifx!csname _@#1@##1!endcsname!relax%
!write16{#1 ##1 not hitherto defined - rerun saving references}%
!changedreferencestrue%
!else%
!expandafter!ifx!csname _@#1@##1!endcsname!temptoken%
!else%
!write16{#1 ##1 reference has changed - rerun saving references}%
!changedreferencestrue%
!fi%
!fi%
!else%
!expandafter!edef!csname _@#1@##1!endcsname{!temptoken}%
!edef!textoutput{!write!references{\global\def\_@#1@##1{!temptoken}}}%
!textoutput%
!fi}}%
!global!def!makecounter#1{!_@@internal@@makelabel{#1}!_@@internal@@makeref{#1}}%
!unsetcatcodes%
}
%
%
%
%
%
\def\turnintolatin#1{\ifcase #1 _\or i\or ii\or iii\or iv\or v\or vi\or vii\or viii\or ix\or x\or xi\or xii\or xiii\or xiv\or xv\or xvi\or xvii\or xviii\or xix\or xx\or xxi\or xxii\or xxiii\or xxiv\or xxv\or xxvi\fi}%
\def\alphanum#1{\ifcase #1 _\or A\or B\or C\or D\or E\or F\or G\or H\or I\or J\or K\or L\or M\or N\or O\or P\or Q\or R\or S\or T\or U\or V\or W\or X\or Y\or Z\fi}%
\newwrite\references%
\ifloadreferences{\catcode`\@=11 \catcode`\_=11 \global\def\_@citation@AnderssonBarbotBeguinZeghib{1}
\global\def\_@citation@Apanasov{2}
\global\def\_@citation@Barbot{3}
\global\def\_@citation@BarbotBonsanteSchlenker{4}
\global\def\_@citation@BarbotBeguinZeghibI{5}
\global\def\_@citation@SanchezI{6}
\global\def\_@citation@SanchezII{7}
\global\def\_@citation@Bonsante{8}
\global\def\_@citation@BonsanteFillastre{9}
\global\def\_@citation@BonsanteMondelloSchlenkerI{10}
\global\def\_@citation@BonsanteMondelloSchlenkerII{11}
\global\def\_@citation@CaffNirSprV{12}
\global\def\_@citation@Carlip{13}
\global\def\_@citation@FillastreSmith{14}
\global\def\_@citation@GilbTrud{15}
\global\def\_@citation@Goldman{16}
\global\def\_@citation@JohnsonMillson{17}
\global\def\_@citation@Kapovich{18}
\global\def\_@citation@KapovichII{19}
\global\def\_@citation@Kourouniotis{20}
\global\def\_@citation@KulkPink{21}
\global\def\_@citation@LabA{22}
\global\def\_@citation@Mess{23}
\global\def\_@citation@Scannell{24}
\global\def\_@citation@SmiRosDT{25}
\global\def\_@citation@SmiFCS{26}
\global\def\_@citation@SmiSLC{27}
\global\def\_@citation@SpruckXiao{28}
\global\def\_@head@Introduction{1}
\global\def\_@subhead@Introduction{1.1}
\global\def\_@proc@MainTheorem{1.1.1}
\global\def\_@head@SpecialLegendrianGeometry{2}
\global\def\_@subhead@RealSpecialLagrangianStructures{2.1}
\global\def\_@proc@WhenTheHolomorphicFormYieldsUnity{2.1.2}
\global\def\_@subhead@TheRealSpecialLegendrianStructure{2.2}
\global\def\_@eqn@RealAndImaginaryComponents{\relax \unhbox \voidb@x \hbox {{\relax \tenrm (2.1)}}}
\global\def\_@eqn@CurvatureOfContactBundle{\relax \unhbox \voidb@x \hbox {{\relax \tenrm (2.2)}}}
\global\def\_@subhead@SpecialLegendrianSubmanifolds{2.3}
\global\def\_@eqn@CodazziMainardi{\relax \unhbox \voidb@x \hbox {{\relax \tenrm (2.3)}}}
\global\def\_@eqn@GraphicalSLCondition{\relax \unhbox \voidb@x \hbox {{\relax \tenrm (2.4)}}}
\global\def\_@proc@Minimal{2.3.2}
\global\def\_@subhead@PositiveSpecialLegendrianSubmanifolds{2.4}
\global\def\_@proc@Dictionary{2.4.1}
\global\def\_@proc@CompactnessPSLC{2.4.2}
\global\def\_@subhead@DegenerateSubmanifolds{2.5}
\global\def\_@proc@OrthogonalToMinimizers{2.5.1}
\global\def\_@eqn@LaplacianOfPhi{\relax \unhbox \voidb@x \hbox {{\relax \tenrm (2.5)}}}
\global\def\_@proc@FormulaForDeltaPhi{2.5.2}
\global\def\_@eqn@FormulaForDeltaXX{\relax \unhbox \voidb@x \hbox {{\relax \tenrm (2.6)}}}
\global\def\_@subhead@RefiningTheSpecialLegendrianCondition{2.6}
\global\def\_@proc@MJIsPositive{2.6.1}
\global\def\_@proc@TotalDegeneracy{2.6.2}
\global\def\_@subhead@TheGeometryOfCurtainSubmanifolds{2.7}
\global\def\_@proc@DescribingCurtainSubmanifolds{2.7.1}
\global\def\_@proc@CurtainSubmanifolds{2.7.2}
\global\def\_@head@HypersurfacesInMinkowskiSpace{3}
\global\def\_@subhead@QuasifuchsianRepresentations{3.1}
\global\def\_@subhead@Stability{3.2}
\global\def\_@proc@JacobiOperator{3.2.1}
\global\def\_@proc@Openness{3.2.3}
\global\def\_@subhead@Compactness{3.3}
\global\def\_@proc@Position{3.3.1}
\global\def\_@proc@UniformlySpacelike{3.3.2}
\global\def\_@proc@CompactnessForHypersurfaces{3.3.3}
\global\def\_@subhead@ExistenceAndUniquenss{3.4}
\global\def\_@proc@MainTheoremWithFoliations{3.4.1}
\global\def\_@head@Bibliography{4}
 }%
\else{\openout\references=references.tex }%
\fi%
%
%
\newcount\headno%
\global\headno=0%
\def\headinfo{\ifinappendices\alphanum\headno\else\the\headno\fi}%
\def\nextheadno{\global\advance\headno by 1 \global\subheadno=0 \global\procno=0 \global\eqnno=0 \headinfo}%
\makecounter{head}%
%
%
\newcount\subheadno%
\global\subheadno=0%
\def\subheadinfo{\headinfo.\the\subheadno}%
\def\nextsubheadno{\global\advance\subheadno by 1 \global\procno=0 \subheadinfo}%
\makecounter{subhead}%
%
%
\newcount\procno%
\global\procno=0%
\def\procinfo{\subheadinfo.\the\procno}%
\def\nextprocno{\global\advance\procno by 1 \procinfo}%
\makecounter{proc}%
%
%
\newcount\figno%
\global\figno=0%
\def\figinfo{\subheadinfo.\the\figno}%
\def\nextfigno{\global\advance\figno by 1 \figinfo}%
\makecounter{fig}%
%
%
\newcount\eqnno%
\global\eqnno=0%
\def\eqninfo{\text{{\rm (\headinfo.\the\eqnno)}}}%
\def\nexteqnno[#1]{\global\advance\eqnno by 1 \eqninfo\hbox{\eqnlabel{#1}}}%
\makecounter{eqn}%
%
%
%
%
%
\def\gobbleeight#1#2#3#4#5#6#7#8{}%
\newcount\citationno%
\global\citationno=0%
\def\citationinfo{\the\citationno}%
\makecounter{citation}%
\newwrite\biblio%
\def\newref#1#2{%
\def\temptext{#2}%
\edef\bibliotextoutput{\expandafter\gobbleeight\meaning\temptext}%
\global\advance\citationno by 1\citationlabel{#1}%
\ifmakebiblio%
    \edef\fileoutput{\write\biblio{\noindent\hbox to 0pt{\hss$[\the\citationno]$}\hskip 0.2em\bibliotextoutput\medskip}}%
    \fileoutput%
\fi}%
\def\cite#1{%
$[\citationref{#1}]$%
\ifmakebiblio%
    \edef\fileoutput{\write\biblio{#1}}%
    \fileoutput%
\fi%
}%
%
%
%
%
\let\mypar=\par%
\edef\Pagetitle={Blank}\headline={\hfil\Pagetitle\hfil}%
\edef\Pagefooter={Blank}\footline={\hfil\Pagefooter\hfil}%
%
%
\newcount\showpagenumflag%
\global\showpagenumflag=0 %
\def\nextoddpage%
{\newpage\ifodd\pageno%
\else\global\showpagenumflag=0 %
\null\vfil\eject%
\global\showpagenumflag=1 %
\fi}%
%
%
\font\headfont=cmb12%
\def\newhead#1[#2]%
{\ifhmode\mypar\fi%
\ifnum\headno=0 \else\goodbreak\bigskip\fi%
{\headfont\noindent\nextheadno\ - #1.}\headlabel{#2}%
\nobreak\medskip}%
%
%
\def\newsubhead#1[#2]%
{\ifhmode\mypar\fi%
\ifnum\subheadno=0 \else\goodbreak\medskip\fi%
{\bf\noindent\nextsubheadno\ - #1.\ }\subheadlabel{#2}}%
%
%
\newif\ifinproclaim%
\global\inproclaimfalse%
\def\proclaim#1{%
\goodbreak\medskip
\bgroup\inproclaimtrue%
\noindent{\bf #1}%
\nobreak\medskip\sl}%
\def\noskipproclaim#1{%
\goodbreak\medskip%
\bgroup\inproclaimtrue%
\noindent{\bf #1}\nobreak\sl}%
\def\endproclaim{\mypar\egroup\nobreak\medskip\ignorespaces}%
%
%
%
\newcount\xpos\newcount\ypos
\def\makelabelgrid{%
\xpos=-5 \ypos=-5 %
\loop\ifnum\xpos<6 %
{\loop\ifnum\ypos<6 %
\def\labeltext{x}%
\ifnum\xpos=0\def\labeltext{+}\fi%
\ifnum\ypos=0\def\labeltext{+}\fi%
\placelabel[\xpos][\ypos]{\labeltext}%
\advance\ypos by 1 %
\repeat}%
\advance\xpos by 1 %
\repeat}%
\def\placelabel[#1][#2]#3{{%
\setbox10=\hbox{\raise #2cm \hbox{\hskip #1cm #3}}%
\ht10=0pt \dp10=0pt \wd10=0pt \box10}}%
%
%
%
%
\def\myitem#1{\noindent\hbox to .5cm{\hfill#1\hss}}%
%
%
%
%
%
%
%
%
%
\font\sansseriften=cmss10%
\font\sansserifseven=cmss7%
\font\sansseriffive=cmss5%
\newfam\sansseriffam%
\textfont\sansseriffam=\sansseriften%
\scriptfont\sansseriffam=\sansserifseven%
\scriptscriptfont\sansseriffam=\sansseriffive%
\def\mathsf{\fam\sansseriffam}%
%
%
%
\font\boldten=cmb10%
\font\boldseven=cmb7%
\font\boldfive=cmb5%
\newfam\mathboldfam%
\textfont\mathboldfam=\boldten%
\scriptfont\mathboldfam=\boldseven%
\scriptscriptfont\mathboldfam=\boldfive%
\def\mathbf{\fam\mathboldfam}%
%
%
%
\font\mycmmiten=cmmi10%
\font\mycmmiseven=cmmi7%
\font\mycmmifive=cmmi5%
\newfam\mycmmifam%
\textfont\mycmmifam=\mycmmiten%
\scriptfont\mycmmifam=\mycmmiseven%
\scriptscriptfont\mycmmifam=\mycmmifive%
\def\hexa#1{\ifcase #1 0\or 1\or 2\or 3\or 4\or 5\or 6\or 7\or 8\or 9\or A\or B\or C\or D\or E\or F\fi}%
\mathchardef\mathi="7\hexa\mycmmifam7B%
\mathchardef\mathj="7\hexa\mycmmifam7C%
%
%
\font\mymsbmten=msbm10 at 8pt%
\font\mymsbmseven=msbm7 at 5.6pt
\font\mymsbmfive=msbm5 at 4pt%
\newfam\mymsbmfam%
\textfont\mymsbmfam=\mymsbmten%
\scriptfont\mymsbmfam=\mymsbmseven%
\scriptscriptfont\mymsbmfam=\mymsbmfive%
\mathchardef\mybeth="7\hexa\mymsbmfam69%
\mathchardef\mygimmel="7\hexa\mymsbmfam6A%
\mathchardef\mydaleth="7\hexa\mymsbmfam6B%
%
%
%
%
\def\proof{{\noindent\bf Proof:\ }}%
\def\remark{{\noindent\bf Remark:\ }}%
\def\qed{~$\square$}%
\def\makeop#1{\global\expandafter\def\csname op#1\endcsname{{\text{#1}}}}%
\def\makeopsmall#1{\global\expandafter\def\csname op#1\endcsname{{\text{\lowercase{#1}}}}}%
%
%
%
\def\minter{\mathop{\cap}}%
%
%
\makeop{Ext}%
\makeop{Int}%
\makeop{Dist}%
\makeop{Diam}%
\makeop{Length}%
%
%
%
%
%
%
%
%
%
\def\minf{\mathop{{\text{Inf}}}}%
%
%
\makeop{Dim}%
\makeop{Ker}%
\makeop{Coker}%
\makeop{Tr}%
\makeop{Adj}%
\makeop{Det}%
\makeop{End}%
\makeop{Lin}%
\makeop{Symm}%
\makeop{Mult}%
%
%
\makeop{dx}%
\makeop{dy}%
\makeop{dz}%
\makeop{dt}%
\makeop{dVol}%
\makeop{dArea}%
\makeop{Supp}%
\makeop{Hess}%
\makeop{Lip}%
%
%
\makeop{Re}%
\makeop{Im}%
\makeop{Arg}%
\makeop{Log}%
\makeop{Exp}%
%
%
\makeopsmall{Cos}%
\makeopsmall{Sin}%
\makeopsmall{Tan}%
\makeopsmall{Sec}%
\makeopsmall{Cosec}%
\makeopsmall{Cot}%
\makeopsmall{ArcCos}%
\makeopsmall{ArcSin}%
\makeopsmall{ArcTan}%
\makeopsmall{ArcSec}%
\makeopsmall{ArcCosec}%
\makeopsmall{ArcCot}%
%
%
\makeopsmall{Cosh}%
\makeopsmall{Sinh}%
\makeopsmall{Tanh}%
\makeopsmall{ArcCosh}%
\makeopsmall{ArcSinh}%
\makeopsmall{ArcTanh}%
%
%
\makeop{Vol}%
\makeop{Area}%
\makeop{Riem}%
\makeop{Ric}%
\makeop{Scal}%
\makeop{Euc}%
\makeop{Imm}%
\makeop{Emb}%
%
%
\makeop{Id}%
\makeop{Ad}%
\makeop{O}%
\makeop{SO}%
\makeop{SL}%
\makeop{GL}%
\makeop{Conf}%
\makeop{Homeo}%
\makeop{Diff}%
\makeop{Isom}%
%
%
\makeop{Ind}%
\makeop{Sig}%
\makeop{Spec}%
%
%
\makeop{Conv}%
\makeop{Max}%
\makeop{Min}%
\makeop{Mod}%
\makeop{Deg}%
\makeop{loc}%
%
%
%
%
%
%
%
%
%
%
%
%
%
 %
%
%
%
%
%
\def\Pagetitle{\hfil}
\def\Pagefooter{\hfil{\myfontdefault\folio}\hfil}
\null \vfill
\def\centre{\rightskip=0pt plus 1fil \leftskip=0pt plus 1fil \spaceskip=.3333em \xspaceskip=.5em \parfillskip=0em \parindent=0em}%
\def\textmonth#1{\ifcase#1\or January\or Febuary\or March\or April\or May\or June\or July\or August\or September\or October\or November\or December\fi}
\font\abstracttitlefont=cmr10 at 14pt {\abstracttitlefont\centre
Constant scalar curvature hypersurfaces\break in $(3+1)$-dimensional GHMC Minkowski spacetimes.\par}
\bigskip
{\centre 26th March 2017\par}
\bigskip
{\centre Graham Smith\par}
\medskip
{\centre Instituto de Matem\'atica,\par
UFRJ, Av. Athos da Silveira Ramos 149,\par
Centro de Tecnologia - Bloco C,\par
Cidade Universit\'aria - Ilha do Fund\~ao,\par
Caixa Postal 68530, 21941-909,\par
Rio de Janeiro, RJ - BRASIL\par}
\bigskip
\noindent{\bf Abstract:\ }We prove that every $(3+1)$-dimensional flat GHMC Minkowski spacetime which is not a translation spacetime or a Misner spacetime carries a unique foliation by spacelike hypersurfaces of constant scalar curvature. In otherwords, we prove that every such spacetime carries a unique time function with isochrones of constant scalar curvature. Furthermore, this time function is a smooth submersion.
%
%
\par
\vfill
\eject
%
\myfontdefault
\def\Pagetitle{\hfil {\rm Constant scalar curvature hypersurfaces...}\hfil}
\def\Pagefooter{\hfil{\myfontdefault\folio}\hfil}
\catcode`\@=11
\def\triplealign#1{\null\,\vcenter{\openup1\jot \m@th %
\ialign{\strut\hfil$\displaystyle{##}\quad$&\hfil$\displaystyle{{}##}$&$\displaystyle{{}##}$\hfil\crcr#1\crcr}}\,}
\def\multiline#1{\null\,\vcenter{\openup1\jot \m@th %
\ialign{\strut$\displaystyle{##}$\hfil&$\displaystyle{{}##}$\hfil\crcr#1\crcr}}\,}
\catcode`\@=12
\makeop{J}%
\makeop{A}%
\makeop{Arctan}%
\makeop{U}%
\makeop{R}%
\makeop{I}%
\makeop{II}%
\makeop{PSO}%
\makeop{inv}%
\makeop{C}%
\newref{AnderssonBarbotBeguinZeghib}{Andersson L., Barbot T., B\'eguin F., Zeghib A., Cosmological time versus CMC time in spacetimes of constant curvature, {\sl Asian Journal of Mathematics}, {\bf 16}, (2012), no. 1, 37--88}
\newref{Apanasov}{Apanasov B. N., Deformations of conformal structures on hyperbolic manifolds, {\sl J. Differential Geom.}, {\bf 35}, (1992), no. 1, 1--20}
\newref{Barbot}{Barbot T., Flat globally hyperbolic spacetimes,{\sl J. Geom. Phys.}, {\bf 53}, (2005), no.2, 123--165}
\newref{BarbotBonsanteSchlenker}{Barbot T., Bonsante F., Schlenker J.-M., Collisions of particles in locally AdS spacetimes I. Local description and global examples, {\sl Comm. Math. Phys.}, {\bf 308}, (2011), no. 1, 147--200}
\newref{BarbotBeguinZeghibI}{Barot T., B\'eguin F., Zeghib A., Prescribing Gauss curvature of surfaces in 3-dimensional spacetimes, Application to the Minkowski problem in the Minkowski space, {\sl Ann. Instit. Fourier.}, {\bf 61}, (2011), no. 2, 511--591}
\newref{SanchezI}{Bernal A. N., S\'anchez M., On smooth Cauchy hypersurfaces and Geroch’s splitting theorem, {\sl Comm. Math. Phys.}, {\bf 243}, (2003), no. 3, 461--470}
\newref{SanchezII}{Bernal A. N., S\'anchez M., Globally hyperbolic spacetimes can be defined as ‘causal’ instead of ‘strongly causal’, {\sl Classical Quantum Gravity}, {\bf 24}, (2007), no. 3, 745--749}
\newref{Bonsante}{Bonsante F., Flat spacetimes with compact hyperbolic Cauchy surfaces, {\sl J. Differential Geom.}, {\bf 69}, (2005), no. 3, 441--521}
\newref{BonsanteFillastre}{Bonsante F., Fillastre F., The equivariant Minkowski problem in Minkowski space, to appear in {\sl Ann. Inst. Fourier}}
\newref{BonsanteMondelloSchlenkerI}{Bonsante F., Mondello G., Schlenker J.-M., A cyclic extension of the earthquake flow, {\sl Geometry \& Topology}, {\bf 17}, (2013), 157--234}
\newref{BonsanteMondelloSchlenkerII}{Bonsante F., Mondello G., Schlenker J.-M., A cyclic extension of the earthquake flow II, arXiv:1208.1738}
\newref{CaffNirSprV}{Caffarelli L., Nirenberg L., Spruck J., Nonlinear second-order elliptic equations. V. The Dirichlet problem for Weingarten hypersurfaces, {\sl Comm. Pure Appl. Math.}, {\bf 41}, (1988), no. 1, 47--70}
\newref{Carlip}{Carlip S., {\sl Quantum gravity in $2+1$ dimensions}, Cambridge Monographs of Mathematical Physics, Cambridge University Press, Cambridge, (1998)}
\newref{FillastreSmith}{Fillastre F., Smith G., Group actions and scattering problems in Teichmueller theory, arXiv:1605.04563}
\newref{GilbTrud}{Gilbarg D., Trudinger N. S., {\sl Elliptic partical differential equations of second order}, Die Grundlehren der mathemathischen Wissenschaften, {\bf 224}, Springer-Verlag, Berlin, New York (1977)}
\newref{Goldman}{Goldman W. M., {\sl Discontinuous groups and the Euler class}, PhD Thesis, University of California, Berkeley, 1980, 138 pp.}
\newref{JohnsonMillson}{Johnson D., Millson J. J., Deformation spaces associated to compact hyperbolic manifolds, {\sl Discrete groups in geometry and analysis (New Haven, Conn., 1984)}, Progr. Math., vol. {\bf 67}, Birkh\"aauser Boston, Boston, MA, 1987, pp. 48--106}
\newref{Kapovich}{Kapovich M., Deformations of representations of discrete subgroups of $\opSO(3,1)$, {\sl Math. Ann.}, {\bf 299}, (1994), no. 2, 341--354}
\newref{KapovichII}{Kapovich M., Kleinian Groups in Higher Dimensions, in {\sl Geometry and Dynamics of Groups and Spaces}, Progress in Mathematics, {bf 265}, 487--564}
\newref{Kourouniotis}{Kourouniotis C., Deformations of hyperbolic structures, {\sl Math. Proc. Cambridge Philos. Soc.}, {\bf 98}, (1985), no. 2, 247--261}
\newref{KulkPink}{Kulkarni R.S., Pinkall U., A canonical metric for M\"obius structures and its applications, {\sl Math. Z.}, {\bf 216}, (1994), no.1, 89--129}
\newref{LabA}{Labourie F., Un lemme de Morse pour les surfaces convexes, {\sl Invent. Math.}, {\bf 141}, (2000), no. 2, 239--297}
\newref{Mess}{Mess G., Lorentz spacetimes of constant curvature, {\sl Geom. Dedicata}, {\bf 126}, (2007), 3--45}
\newref{Scannell}{Scannell K. P., Infinitesimal deformations of some $\opSO(3,1)$ lattices, {\sl Pacific J. Math.}, {\bf 194}, (2000), no. 2, 455--464}
\newref{SmiRosDT}{Rosenberg H., Smith G., Degree Theory of Immersed Hypersurfaces, arXiv:1010.1879}
\newref{SmiFCS}{Smith G., Moduli of Flat Conformal Structures of Hyperbolic Type, {\sl Geom. Dedicata}, {\bf 154}, (2011), no. 1, 47--80}
\newref{SmiSLC}{Smith G., Special Lagrangian curvature, Math. Annalen, {\bf 335}, (2013), no. 1, 57--95}
\newref{SpruckXiao}{Spruck J., Xiao L., Convex Spacelike Hypersurfaces of Constant Curvature in de Sitter Space, arXiv:1203.5739}
%
%
%
\newhead{Introduction}[Introduction]
\newsubhead{Introduction}[Introduction]
Let $\Bbb{R}^{d,1}$ denote $(d+1)$-dimensional Minkowski space, that is $\Bbb{R}^{d+1}$ furnished with the semi-riemannian metric
$$
ds^2 := dx_1^2 + ... + dx_d^2 - dx_{d+1}^2.
$$
For the purposes of this paper, a $(d+1)$-dimensional {\sl Minkowski spacetime} is a semi-riemannian manifold $X$ which is everywhere locally isometric to $\Bbb{R}^{d,1}$. Minkowski spacetimes are of interest in cosmology as the simplest possible solutions of Einstein's equations. However, since the pioneering work \cite{Mess} of Mess, they have also found deep and broad applications in the study of Teichm\"uller theory and its higher-dimensional analogous.
\par
In this paper, we will be concerned with time functions defined over certain types of Minkowski spacetimes known as GHMC Minkowski spacetimes (see below). Here, a {\sl time function} is defined to be a real-valued submersion all of whose level sets are spacelike. It is often possible to construct time functions which have interesting geometric properties and these, in turn, can be useful in studying the physics or the geometry of the ambient spacetime. For example, in \cite{AnderssonBarbotBeguinZeghib}, Andersson, Barbot, B\'eguin \& Zeghib study smooth time-functions over constant curvature\footnote{${}^*$}{that is, de-Sitter, Minkowski and anti de-Sitter.} GHMC spacetimes with level sets of constant mean curvature, thus proving the existence over such spacetimes of so-called ``York'' time functions, which are known to be of considerable use in general relativity (c.f. \cite{Carlip}). Likewise, in \cite{BarbotBonsanteSchlenker} and \cite{BarbotBeguinZeghibI}, Barbot, B\'eguin \& Zeghib construct smooth time functions over $(2+1)$-dimensional GHMC Minkowski and anti de-Sitter spacetimes with level sets of constant extrinsic curvature, and in \cite{BonsanteMondelloSchlenkerI} and \cite{BonsanteMondelloSchlenkerII}, Bonsante, Mondello \& Schlenker use these time functions to provide fascinating new insights into the earthquake and grafting maps of classical Teichm\"uller theory. We refer the interested reader to our review \cite{FillastreSmith} for more details of these and other related constructions.
\par
In contrast to the above mentioned results, in \cite{BonsanteFillastre}, Bonsante \& Fillastre show that smooth time functions with level sets of constant extrinsic curvature do not always exist for higher dimensional ambient spacetimes. This can be understood as a manifestation of the subtle manner in which constant extrinsic curvature is actually a uniformly elliptic condition for surfaces, but is no longer so in higher dimensions. For this reason, in \cite{SmiFCS}, we introduced the so-called ``special lagrangian curvature'' which, by involving the Calabi-Yau structure of the tangent bundle of the ambient space, yields an alternative generalisation of extrinsic curvature which continues to be uniformly elliptic even in higher dimensions. In particular, when the ambient spacetime is $(3+1)$-dimensional, the ``special lagrangian curvature'' of any spacelike hypersurface is none other than its {\sl scalar curvature} which, we recall, is defined over a riemannian manifold $M$ by
$$
S:=\frac{-1}{m(m-1)}g^{ik}R_{ijk}{}^j,
$$
where $m$ here denotes the dimension of $M$, $g$ denotes its metric, $R$ denotes its Riemann curvature tensor, and the summation convention is implied.\footnote{${}^\dag$}{The normalisation is chosen here so that the scalar curvature of the unit sphere in Euclidean space is equal to $1$.} Using this curvature notion, in \cite{SmiFCS}, we partially complemented the work \cite{AnderssonBarbotBeguinZeghib} of Andersson, Barbot, B\'eguin \& Zeghib by constructing smooth time functions over open subsets of GHMC de-Sitter spacetimes with level sets of constant scalar curvature. In the present work, we extend this result to the Minkowski case, that is, we construct smooth time functions over $(3+1)$-dimensional GHMC Minkowski spacetimes with level sets of constant scalar curvature.
\par
Before stating the main result of this paper, it is necessary to properly introduce the class of GHMC Minkowski spacetimes. Although this class arises from fairly natural physical considerations, it is nonetheless rather unfamiliar to most geometers, and its definition therefore requires a brief detour. We begin with the concept of causality. A tangent vector of $X$ is said to be {\sl spacelike}, {\sl timelike} or {\sl null} according to whether its norm-squared is positive, negative or null respectively, and is said to be {\sl causal} whenever it is either timelike or null. A continuously differentiable, embedded curve in $X$ is said to be {\sl causal} whenever all of its tangent vectors are causal. Finally, the spacetime $X$ is itself said to be {\sl causal} whenever it contains no closed causal curve. This condition reflects the physical requirement that it is not possible to move into the past by travelling into the future.
\par
We henceforth suppose that $X$ is oriented, time-oriented and causal. The {\sl past} of any given point $x$ of $X$ is then defined to be the closed set of all initial points of future-oriented, causal curves terminating in $x$. Likewise, the {\sl future} of that point is defined to be the closed set of all terminal points of future-oriented, causal curves starting at $x$. The spacetime $X$ is said to be {\sl globally hyperbolic} whenever the intersection of the past of any point with the future of any other point is compact (c.f. \cite{SanchezI} and \cite{SanchezII}).
\par
A {\sl Cauchy hypersurface} in $X$ is a spacelike hypersurface $\Sigma$ which intersects every inextensible causal curve exactly once. In \cite{SanchezI} and \cite{SanchezII}, Bernal \& Sanchez show that an oriented and time-oriented spacetime is globally hyperbolic if and only if it contains a smooth Cauchy hypersurface. Although the Cauchy hypersurface is trivially not unique, all smooth Cauchy hypersurfaces of a given globally hyperbolic spacetime are diffeomorphic to one another, and the spacetime itself is diffeomorphic to the Cartesian product of any such hypersurface with $\Bbb{R}$. A globally hyperbolic spacetime is said to be {\sl Cauchy compact} whenever its Cauchy hypersurface is compact.
\par
Finally, a globally hyperbolic spacetime $X$ is said to be {\sl maximal} whenever there exists no isometric embedding $e$ of $X$ into a strictly larger spacetime $\tilde{X}$ such that the image under $e$ of a Cauchy hypersurface in $X$ is also a Cauchy hypersurface in $\tilde{X}$. In order to understand this property, consider the future cone in $\Bbb{R}^{d,1}$,
$$
\opC^{d,1} := \left\{ x\in\Bbb{R}^{d,1}\ |\ \|x\|^2<0,\ x_{d+1}>0\right\}.
$$
This space is globally hyperbolic and maximal, even though it isometrically embeds into the strictly larger space $\Bbb{R}^{d,1}$. Indeed, its Cauchy hypersurface is given by
$$
\Sigma^d := \left\{ x\in\Bbb{R}^{d,1}\ |\ \|x\|^2=-1,\ x_{d+1}>0\right\},
$$
which is not a Cauchy hypersurface of $\Bbb{R}^{d,1}$.
\par
We say that an oriented and time oriented spacetime is {\sl GHMC} whenever it is Globally Hyperbolic, Maximal and Cauchy-compact. Building on the work \cite{Mess} of Mess, a classification of all GHMC Minkowski spacetimes was initiated by Bonsante in \cite{Bonsante} and completed by Barbot in \cite{Barbot}. Within this classification, two exceptional families stand out, namely the {\sl translation spacetimes} and the {\sl Misner spacetimes}. A GHMC translation spacetime $X$ is one whose universal cover is the whole of $\Bbb{R}^{d,1}$. Up to a finite cover, such a spacetime has the form
$$
X = T^d\times\Bbb{R},
$$
furnished with the metric
$$
g = dx_1^2 + ... + dx_d^2 - dx_{d+1}^2,
$$
where $T^d$ is the quotient of $\Bbb{R}^d$ by some cocompact lattice $\Lambda$, and it carries the foliation $(T^d\times\left\{t\right\})_{t\in\Bbb{R}}$, all of whose leaves have vanishing scalar curvature. A GHMC Misner spacetime $Y$ is, up to reversal of the temporal orientation, one whose universal cover is $\Bbb{R}^{d-1}\times C^{1,1}$, where $C^{1,1}$ is the future cone in $\Bbb{R}^{1,1}$. Up to a finite cover, every such spacetime has the form
$$
Y = T^d\times]0,\infty[,
$$
furnished with the metric
$$
g = dx_1^2 + ... + dx_{d-1}^2 + t^2dx_d^2 - dx_{d+1}^2,
$$
where $T^d$ is again the quotient of $\Bbb{R}^d$ by some cocompact lattice $\Lambda$. However, in this case, it is worth observing that the projection of the $x_d$-axis need not necessarily yield a closed curve in $T^d$. Such a spacetime carries the foliation $(T^d\times\left\{t\right\})_{t>0}$, all of whose leaves also have vanishing scalar curvature.
\par
In both of the above cases, it follows from the geometric maximum principle that there exist no immersed spacelike hypersurfaces of constant non-zero scalar curvature. For all other GHMC Minkowski spacetimes, we have
\proclaim{Theorem \nextprocno}
\noindent Let $X$ be a $(3+1)$-dimensional GHMC Minkowski spacetime which is not a translation spacetime or a Misner spacetime. There exists a unique smooth submersion $T:X\rightarrow]0,\infty[$ such that, for all $t\in]0,\infty[$, the level set $T^{-1}(\left\{t\right\})$ is a convex Cauchy hypersurface of constant scalar curvature equal to $(-t^2)$.
\endproclaim
\proclabel{MainTheorem}
\remark By the classification \cite{Barbot} of Barbot, every GHMC Minkowski spacetime which is not a translation spacetime or a Misner spacetime is, up to reversal of the temporal orientation and up to a finite cover, a twisted product of what we choose to call a ``Bonsante spacetime'' with a Euclidean torus (c.f. Section \subheadref{QuasifuchsianRepresentations}, below). Theorem \procref{MainTheorem} is then a straightforward consequence of the corresponding result for Bonsante spacetimes, and thus follows immediately from Theorem \procref{MainTheoremWithFoliations}, below.
\medskip
\remark In fact, an analogous results holds in all dimensions for submersions with level hypersurfaces of constant ``special lagrangian curvature''. In particular, when the ambient space is $(2+1)$-dimensional, we recover the result \cite{BarbotBeguinZeghibI} of Barbot, B\'eguin \& Zeghib. However, in higher dimensions, the ``special lagrangian curvature'' no longer has an elementary expression in terms of classical curvature notions and, for this reason, we shall not discuss this further here.
\medskip
The author is grateful to Fran\c{c}ois Fillastre for long and enlightening conversations on the subject of semi-riemannian geometry, without which this work would not have been realized. The author is also grateful to Thierry Barbot for helpful comments concerning earlier drafts of this paper. This paper was written as part of the project MATH AMSUD 2017, Project No. 38888QB - GDAR.
\newhead{Special legendrian geometry}[SpecialLegendrianGeometry]
\newsubhead{Real special lagrangian structures}[RealSpecialLagrangianStructures]
The proof of Theorem \procref{MainTheorem} is based on a compactness result for families of hypersurfaces of constant scalar curvature. This result in turn depends on the special legendrian structure of the bundle of future-oriented, unit, timelike tangent vectors over $\Bbb{R}^{d,1}$. In order to understand this structure, we first consider the fully integrable case. Thus, denote by $\langle\cdot,\cdot\rangle$ the Minkowski metric over $\Bbb{R}^{d,1}$, that is
$$
\langle X,Y\rangle = X^1Y^1 + ... + X^dY^d - X^{d+1}Y^{d+1},
$$
and consider the following structures defined over the cartesian product $\Bbb{R}^{d,1}\times\Bbb{R}^{d,1}$.
$$\eqalign{
\omega((X_r,X_i),(Y_r,Y_i)) &:= \langle X_r,Y_i\rangle - \langle Y_r,X_i\rangle,\cr
\opJ(X_r,X_i) &:= (-X_i,X_r),\ \text{and}\cr
\opR(X_r,X_i) &:= (X_r,-X_i).\cr}
$$
Observe that $\omega$ is antisymmetric and non-degenerate,
$$\eqalign{
\opJ^2 &= -\opId,\ \text{and}\cr
\opR^2 &= \opId.\cr}
$$
The form $\omega$ is the {\sl symplectic} structure, and the maps $\opJ$ and $\opR$ are, respectively, the {\sl complex} and {\sl real} structures. They are related to one another by
$$\eqalign{
\omega(\opJ\cdot,\opJ\cdot) &= \omega(\cdot,\cdot),\cr
\omega(\opR\cdot,\opR\cdot) &= -\omega(\cdot,\cdot),\ \text{and}\cr
\opJ\opR + \opR\opJ &= 0.\cr}
$$
We call the structure defined by the triple $(\omega,\opJ,\opR)$ a {\sl real special lagrangian structure}. Its symmetry group is given by the action of $\opO(d,1)$ on $\Bbb{R}^{d,1}\times\Bbb{R}^{d,1}$ which, for any matrix $M$, maps the pair $(X_r,X_i)$ to the pair $(MX_r,MX_i)$.
\par
The real special lagrangian structure yields various auxiliary structures which play a key role in the sequel. First, there is a semi-riemannian metric of signature $(2d,2)$, given by
$$
g((X_r,X_i),(Y_r,Y_i)) := \omega((X_r,X_i),\opJ(Y_r,Y_i)) =
\langle X_r,Y_r\rangle + \langle X_i,Y_i\rangle.
$$
Next, there is another semi-riemannian metric, this time of signature $(d+1,d+1)$, given by
$$
m((X_r,X_i),(Y_r,Y_i)) := -\omega((X_r,X_i),\opR(Y_r,Y_i)) = \langle X_r,Y_i\rangle + \langle Y_r,X_i\rangle.
$$
In particular, over every lagrangian subspace of $\Bbb{R}^{d,1}\times\Bbb{R}^{d,1}$, $m$ satisfies
$$
m((X_r,X_i),(Y_r,Y_i)) := 2\langle X_r,Y_i\rangle = 2\langle Y_r,X_i\rangle.
$$
Finally, there is, up to a choice of sign, a unique complex $(d+1)$-form (complex with respect to the structure $\opJ$) whose restriction to the real subspace $\Bbb{R}^{d,1}\times\left\{0\right\}$ coincides with the volume form of the metric $g$. This form, which we denote by $\hat{\Omega}$ is also described explicitly as follows. Let $e_1,...,e_{d+1}$ be an orthonormal basis\footnote*{Since the metric of $\Bbb{R}^{d,1}$ has signature $(d,1)$, we take this to mean that $e_1$,...,$e_d$ are spatial and $e_{d+1}$ is temporal.} of $\Bbb{R}^{d,1}$ and, for all $k$, denote $f_k:=\opJ e_k$. Let $e^1,...,e^{d+1},f^1,...,f^{d+1}$ be the basis dual to $e_1,...,e_{d+1},f_1,...,f_{d+1}$ and, for all $k$, denote
$$
dz^k := e^k + if^k.
$$
It follows from the definitions that
$$
\hat{\Omega} = \pm dz^1\wedge...\wedge dz^{d+1}.
$$
In particular, $\hat{\Omega}$ is independent of the orthonormal basis of $\Bbb{R}^{d,1}$ chosen.
\proclaim{Lemma \nextprocno}
\myitem{(1)} The restriction of $g$ to any lagrangian subspace of $\Bbb{R}^{d,1}\times\Bbb{R}^{d,1}$ has signature $(d,1)$.
\medskip
\myitem{(2)} If $e_1,...,e_{d+1}$ is an orthonormal basis of some lagrangian subspace of $\Bbb{R}^{d,1}\times\Bbb{R}^{d,1}$, then
$$
\left|\hat{\Omega}(e_1,...,e_{d+1})\right| = 1.
$$
\endproclaim
\proof Let $E$ be a lagrangian subspace of $\Bbb{R}^{d,1}\times\Bbb{R}^{d,1}$ and denote $F:=\opJ E$. Observe that $E$ and $F$ are mutually orthogonal with respect to $g$. It follows that $\Bbb{R}^{d,1}\times\Bbb{R}^{d,1}=E\oplus F$ and that the restrictions of $g$ to $E$ and $F$ are non-degenerate. Furthermore, since $\opJ$ restricts to an isometry from $E$ to $F$ with respect to $g$, the restrictions of $g$ to these two subspaces both have the same signature. In particular, if this signature is equal to $(p,q)$ then, by orthogonality of $E$ and $F$ again, $(2p,2q)=(2d,2)$, and $(1)$ follows.
\par
Now let $A$ be any isometry sending $E$ into the real subspace $\Bbb{R}^{d,1}\times\left\{0\right\}$. This map extends to a unique unitary map of $\Bbb{R}^{d,1}\times\Bbb{R}^{d,1}$ to itself. In particular,
$$
\left|\hat{\Omega}(e_1,...,e_{d+1})\right|=\left|\hat{\Omega}(Ae_1,...,Ae_{d+1})\right|.
$$
However, if $e_1,...,e_{d+1}$ is an orthonormal basis of $E$, then $Ae_1,...,Ae_{d+1}$ is an orthonormal basis of $\Bbb{R}^{d,1}\times\left\{0\right\}$ so that, by definition,
$$
\hat{\Omega}(Ae_1,...,Ae_{d+1}) = \pm 1,
$$
and $(2)$ follows.\qed
\proclaim{Corollary \nextprocno}
\noindent If $L$ is a lagrangian subspace over which
$$
\opIm\left(e^{i\theta}\hat{\Omega}\right)=0
$$
for some $\theta$ then, for any orthonormal basis $e_1,...,e_{d+1}$ of $L$,
$$
\opRe\left(e^{i\theta}\hat{\Omega}(e_1,...,e_{d+1})\right)=\pm 1.
$$
\endproclaim
\proclabel{WhenTheHolomorphicFormYieldsUnity}
\newsubhead{The real special legendrian structure}[TheRealSpecialLegendrianStructure]
Let $\opU^+\Bbb{R}^{d,1}$ denote the bundle of future-oriented, unit, timelike tangent vectors over $\Bbb{R}^{d,1}$, let $M:=M^d$ denote its total space, that is
$$
M:=M^d:=\opU^+\Bbb{R}^{d,1} = \left\{ (x,y)\ |\ \|y\|^2 = -1,\ y_{d+1}>0\right\}.
$$
The tangent space to $M$ at the point $(x,y)$ is given by
$$
T_{(x,y)}M = \left\{ (X_r,X_i)\ |\ \langle X_i,y\rangle = 0 \right\}.
$$
Observe now that the symplectic form $\omega$ defined in Section \subheadref{RealSpecialLagrangianStructures} is the exterior derivative of the Liouville form
$$
\lambda_{(x,y)}(X_r,X_i) := \langle X_r,y\rangle.
$$
The restriction of this form to $TM$ makes $M$ into a contact manifold. We denote the resulting contact bundle by $\alpha$. Its fibre at the point $(x,y)$ is given by
$$
\alpha_{(x,y)} := \left\{ (X_r,X_i)\ |\ \langle X_r,y\rangle=\langle X_i,y\rangle = 0 \right\}.
$$
\par
Since the maps $\opJ$ and $\opR$ both preserve $\alpha$ they define, together with $\omega:=d\lambda$, a real special lagrangian structure over every fibre of the bundle. We say that the manifold $M$ carries a {\sl real special legendrian structure}. Of particular interest to us will be the real and imaginary subspaces of $\alpha_{(x,y)}$ given by
$$\eqalign{
\Cal{R}_{(x,y)} &:= \left\{X\ |\ X=\opR X\right\} = \left\{(X,0)\ |\ \langle X,y\rangle = 0\right\},\ \text{and}\cr
\Cal{I}_{(x,y)} &:= \left\{X\ |\ X=-\opR X\right\} = \left\{(0,X)\ |\ \langle X,y\rangle = 0\right\}.\cr}
\eqnum{\nexteqnno[RealAndImaginaryComponents]}
$$
From the point of view of the unitary bundle, $\Cal{R}_{(x,y)}$ and $\Cal{I}_{(x,y)}$ are respectively the horizontal and vertical subspaces of $\alpha_{(x,y)}$. Furthermore, each of these spaces naturally identifies with the orthogonal complement of $y$ in $\Bbb{R}^{d,1}$ and, in particular, they each identify with one another.
\par
The forms $g$ and $m$ are defined over each fibre of $\alpha$ as before. Significantly, over each fibre, $g$ now defines a riemannian metric, and $m$ defines a metric of signature $(d,d)$. Let $\Omega$ denote the unique (up to choice of sign) complex $d$-form which coincides with the volume form of the metric $g$ over the subspace
$$
\Cal{R}_{(x,y)} = \left\{ (X,0)\ |\ \langle X,y\rangle = 0 \right\}.
$$
It is straightforward to show that at every point $(x,y)$ of $M$,
$$
\Omega_{(x,y)} = \pm c_{(y,0)}\hat{\Omega},
$$
where $c_{(y,0)}$ here denotes the operator of contraction by the vector $(y,0)$.
\par
Let $D$ denote the restriction to $M$ of the standard differentiation operator over $\Bbb{R}^{d,1}\times\Bbb{R}^{d,1}$. Let $\overline{\nabla}$ denote its composition with orthogonal projection onto the distribution $\alpha$. Let $X:=(X_r,X_i)$ and $Y:=(Y_r,Y_i)$ be tangent vector fields over $M$ taking values in $\alpha$. Recalling that $(y,0)$ and $(0,y)$ are temporal vectors, and therefore have negative norm-squared, we have
$$\eqalign{
\overline{\nabla}_X Y &= D_XY + g\left((y,0),D_XY\right)(y,0) + g\left((0,y),D_XY\right)(0,y)\cr
&=D_XY - g\left(Y,D_X(y,0)\right)(y,0) - g\left(Y,D_X(0,y)\right)(0,y)\cr
&=D_XY - \langle Y_r,X_i\rangle(y,0) - \langle Y_i,X_i\rangle(0,y).\cr}
$$
The {\sl shape operator} of $\alpha$ is defined to be the difference between $D$ and $\overline{\nabla}$, that is
$$
\overline{A}_X Y :=  D_X Y - \overline{\nabla}_X Y = \langle Y_r,X_i\rangle(y,0) + \langle Y_i,X_i\rangle(0,y).
$$
The curvature of the distribution $\alpha$ is then given by
$$\eqalign{
g\left(\overline{R}_{XY}Z,W\right)
&=g\left(\overline{A}_XW,\overline{A}_YZ\right) - g\left(\overline{A}_XZ,\overline{A}_YW\right)\cr
&=-\langle W_r,X_i\rangle\langle Z_r,Y_i\rangle - \langle W_i,X_i\rangle\langle Z_i,Y_i\rangle\cr
&\qquad + \langle Z_r,X_i\rangle\langle W_r,Y_i\rangle + \langle Z_i,X_i\rangle\langle W_i,Y_i\rangle,\cr}
\eqnum{\nexteqnno[CurvatureOfContactBundle]}
$$
where $X$, $Y$, $Z$ and $W$ are vectors in $\alpha$. This formula will play a key role in the computations that follow.
\par
Finally, straightforward computations show that for all vectors $X$ in $\alpha$,
$$
(\overline{\nabla}\opJ)X = (\overline{\nabla}\opR)X = 0,
$$
for all vectors $X$ and $Y$ in $\alpha$,
$$
(\overline{\nabla}g)(X,Y) = (\overline{\nabla}\omega)(X,Y) = (\overline{\nabla}m)(X,Y) = 0,
$$
and for all vectors $X,Y_1,...,Y_d$ in $\alpha$,
$$
(\overline{\nabla}_X\Omega)(Y_1,...,Y_d) = 0.
$$
\newsubhead{Special legendrian submanifolds}[SpecialLegendrianSubmanifolds]
Let $\hat{\Sigma}\subseteq M$ be a $d$-dimensional embedded submanifold. We say that $\hat{\Sigma}$ is {\sl legendrian} whenever its tangent space is always contained in $\alpha$ and
$$
\omega|_{T\hat{\Sigma}} = 0.
$$
Let $\nabla$ be the Levi-Civita covariant derivative of $\hat{\Sigma}$. The shape operator $\hat{A}$ of $\hat{\Sigma}$ is defined by
$$
\hat{A}_X Y := \overline{\nabla}_X Y - \nabla_X Y,
$$
where $X$ and $Y$ are tangent vector fields over $\hat{\Sigma}$. Observe, in particular, that $\hat{A}_X Y$ is always normal to $\hat{\Sigma}$. The second fundamental form of $\hat{\Sigma}$ is then given by
$$
\opII(X,Y,Z) := g(\hat{A}_X Y, \opJ Z) = -\omega(\hat{A}_X Y, Z).
$$
In particular, since $\hat{\Sigma}$ is legendrian, $\opII$ is symmetric in all three variables. Furthermore, we have the following Codazzi-Mainardi equations.
\proclaim{Lemma \nextprocno}
\noindent For all vectors $X$, $Y$, $Z$ and $W$ tangent to $\hat{\Sigma}$,
$$
(\nabla_W\opII)(X,Y,Z) = (\nabla_X\opII)(W,Y,Z) + g(\overline{\opR}_{WX}Y,\opJ Z).\eqnum{\nexteqnno[CodazziMainardi]}
$$
\endproclaim
\proof Indeed, consider a point $p\in\hat{\Sigma}$ and tangent vector fields $X$, $Y$, $Z$ and $W$ which are all parallel at this point. Then
$$\eqalign{
(\nabla_W\opII)(X,Y,Z) &= W\opII(X,Y,Z)\cr
&=-W\omega(\overline{\nabla}_X Y,Z)\cr
&=-\omega(\overline{\nabla}_W\overline{\nabla}_XY,Z) - \omega(\overline{\nabla}_XY,\overline{\nabla}_WZ).\cr}
$$
Since both $Y$ and $Z$ are parallel at $p$, both $\overline{\nabla}_XY$ and $\overline{\nabla}_WZ$ are elements of the lagrangian subspace $\opJ T_p\hat{\Sigma}$, and the second term on the right hand side therefore vanishes. Since both $X$ and $W$ are both parallel at $p$, their commutator $[W,X]$ vanishes at this point, and the result now follows by symmetry.\qed
\medskip
For $\theta\in\Bbb{R}$, we say that the legendrian submanifold $\hat{\Sigma}$ is {\sl $\theta$-special legendrian} whenever
$$
\opIm(e^{i\theta}\Omega)|_{T\hat{\Sigma}} = 0,
$$
and we call $\theta$ the {\sl special legendrian angle} of $\hat{\Sigma}$. Observe, in particular, that if $T\hat{\Sigma}$ is the graph of the linear map $B:\Cal{R}\rightarrow\Cal{I}$, then $\hat{\Sigma}$ is $\theta$-special lagrangian if and only if
$$
\opIm(e^{i\theta}\opDet(\opId + iB)) = 0.\eqnum{\nexteqnno[GraphicalSLCondition]}
$$
We will make good use of this property in the sequel.
\proclaim{Lemma \nextprocno}
\noindent If $\hat{\Sigma}$ is $\theta$-special legendrian for some $\theta$ then, for any orthonormal basis $e_1,...,e_d$ of $T\hat{\Sigma}$,
$$
\sum_{a=1}^d \opII(\cdot,e_a,e_a) = 0.
$$
\endproclaim
\proclabel{Minimal}
\remark Observe that this is not the same as minimality of $\hat{\Sigma}$, since it is still possible for the mean curvature vector of this submanifold to be non-trivial in the direction of the Reeb vector field of $M$.
\medskip
\proof Indeed, consider a point $p$ in $\hat{\Sigma}$. Let $e_1,...,e_d$ be an orthonormal basis of $T_p\hat{\Sigma}$ and extend this to a frame of $\hat{\Sigma}$ which is parallel at $p$. For any other tangent vector $X$ to $\hat{\Sigma}$ at $p$,
$$\eqalign{
X\Omega(e_1,...,e_d) &= (\overline{\nabla}_X\Omega)(e_1,...,e_d) + \sum_{a=1}^d\Omega(e_1,...,\overline{\nabla}_Xe_a,...,e_d)\cr
&=\sum_{a=1}^d\Omega(e_1,...,\hat{A}_Xe_a,...,e_d).\cr}
$$
Since $\hat{A}_Xe_a$ is normal to $T_p\hat{\Sigma}$ for all $a$, it follows that
$$\eqalign{
X\Omega(e_1,...,e_d)
&= \sum_{a,b=1}^d\Omega(e_1,...,g(\hat{A}_X e_a,\opJ e_b)\opJ e_b,...,e_d)\cr
&= i\left(\sum_{a=1}^d\opII(X,e_a,e_a)\right)\Omega(e_1,...,e_d),\cr}
$$
where the second equality follows since $\Omega$ is $\Bbb{C}$-multilinear. Taking the imaginary parts now yields
$$\eqalign{
X\opIm(e^{i\theta}\Omega(e_1,...,e_d))
&=\left(\sum_{a=1}^d\opII(X,e_a,e_a)\right)\opRe(e^{i\theta}\Omega(e_1,...,e_d))\cr
&=\pm\sum_{a=1}^d\opII(X,e_a,e_a),\cr}
$$
where the second equality here follows by Corollary \procref{WhenTheHolomorphicFormYieldsUnity}. The result follows.\qed
\newsubhead{Positive special legendrian submanifolds}[PositiveSpecialLegendrianSubmanifolds]
We say that a special legendrian submanifold $\hat{\Sigma}$ is {\sl positive} whenever
$$
m|_{T\hat{\Sigma}} \geq 0.
$$
Although this may seem like a strong restriction, positive special legendrian submanifolds arise naturally as lifts of certain hypersurfaces of $\Bbb{R}^{d,1}$. As such, they will play a key role in our study of constant scalar curvature hypersurfaces in $(3+1)$-dimensional Minkowski space. Indeed, recall that, by the Gauss-Codazzi equations, the scalar curvature of an embedded hypersurface in $\Bbb{R}^{3,1}$ is given by
$$
S:=-\frac{1}{3}(\lambda_1\lambda_2 + \lambda_1\lambda_3 + \lambda_2\lambda_3),
$$
where $\lambda_1$, $\lambda_2$ and $\lambda_3$ its principal curvatures. We now have
\proclaim{Lemma \nextprocno}
\noindent Let $\Sigma$ be an embedded spacelike hypersurface in $\Bbb{R}^{3,1}$, let $N:\Sigma\rightarrow\Bbb{H}^3$ be its future-oriented, unit, normal vector field and define the embedded hypersurface $\hat{\Sigma}$ by
$$
\hat{\Sigma} := \left\{(x,N(x))\ |\ x\in\Sigma\right\}.
$$
\myitem{(1)} $\hat{\Sigma}$ is legendrian;
\medskip
\myitem{(2)} if $\Sigma$ has constant scalar curvature equal to $(-1/3)$, then $\hat{\Sigma}$ is $\pi/2$-special legendrian;
\medskip
\myitem{(3)} if $\Sigma$ is complete, then so too is $\hat{\Sigma}$; and
\medskip
\myitem{(4)} if $\Sigma$ is locally strictly convex, then $\hat{\Sigma}$ is positive.
\endproclaim
\proclabel{Dictionary}
\remark In fact, an analogous result holds in all dimensions. For example, when $d=2$, special legendrian submanifolds correspond to surfaces of constant extrinsic curvature. However, in higher dimensions, the relevant notion of curvature is more involved, and we refer the reader to \cite{SmiSLC} for details.
\medskip
\proof Consider a point $x$ of $\Sigma$ and denote $y:=N(x)$. Let $\Cal{R}$ and $\Cal{I}$ denote respectively the real and imaginary subspaces of $\alpha_{(x,y)}$ as defined in \eqnref{RealAndImaginaryComponents}. Recall that both $\Cal{R}$ and $\Cal{I}$ identify with the orthogonal complement of $y$ in $\Bbb{R}^{3,1}$, that is, the tangent space to $\Sigma$ at $x$. In particular, via this identification, $T_p\hat{\Sigma}$ itself identifies with the graph of the shape operator of $\Sigma$ at $x$, which we now denote by $A$. For $X$ and $Y$ tangent to $\Sigma$, we now have
$$
\omega((X,AX),(Y,AY)) = \langle X,AY\rangle - \langle AX,Y\rangle = 0,
$$
so that $\hat{\Sigma}$ is legendrian. This proves $(1)$. Next, we have
$$
\opDet(\opId + iA) = 1 + i\opTr(A) - \sigma_2(A) - i\opDet(A),
$$
so that $\opIm(e^{i\pi/2}\opDet(\opId + iA))$ vanishes whenever $\sigma_2(A)=1$. This proves $(2)$. The metric over $\hat{\Sigma}$ is given by
$$
g\left((X,AX),(Y,AY)\right) = \langle X,Y\rangle + \langle AX,AY\rangle,
$$
so that $\hat{\Sigma}$ is complete whenever $\Sigma$ is. This proves $(3)$. Finally,
$$
m((X,AX),(Y,AY)) = 2\langle X,AY\rangle,
$$
so that $m$ is positive definite whenever $A$ is positive definite. This proves $(4)$.\qed
\medskip
Significantly, the hypothesis of positivity makes it straightforward to prove strong compactness results for special legendrian submanifolds. In order to state these results, we first require some terminology. Thus, a {\sl pointed}, embedded submanifold is defined to be a pair $(\hat{\Sigma},p)$ where $\hat{\Sigma}$ is an embedded submanifold and $p$ is a point of $\hat{\Sigma}$. A sequence $(\hat{\Sigma}_m,p_m)$ of complete, pointed, embedded submanifolds is said to {\sl converge} towards the complete, pointed, embedded submanifold $(\hat{\Sigma}_\infty,p_\infty)$ whenever there exists a sequence $(\phi_m)$ of smooth maps from $\hat{\Sigma}_\infty$ into $M$ with the following properties.
\medskip
\myitem{(1)} $\phi_m(p_\infty) = p_m$ for all $m$; and
\medskip
\noindent for every compact subset $K$ of $\hat{\Sigma}_\infty$, there exists $m_0\geq 0$ such that
\medskip
\myitem{(2)} for all $m\geq m_0$, $\phi_m$ defines an embedding over $K$ whose image is contained in $\hat{\Sigma}_m$; and
\medskip
\myitem{(3)} the subsequence $(\phi_m)_{m\geq m_0}$ converges in the $C^\infty$ sense over $K$ to the identity map.
\medskip
\noindent In \cite{SmiSLC}, we prove the following compactness theorem.
\proclaim{Theorem \nextprocno}
\noindent Let $(\theta_m)$ be a sequence of real numbers converging to $\theta_\infty$. For all $m$, let $(\hat{\Sigma}_m,p_m)$ be a complete, pointed, positive, $\theta_m$-special legendrian submanifold of $M$. If there exists a compact subset $K$ of $M$ such that $p_m\in K$ for all $m$, then there exists a complete, pointed, positive, $\theta_\infty$-special legendrian submanifold $(\hat{\Sigma}_\infty,p_\infty)$ towards which the sequence $(\hat{\Sigma}_m,p_m)$ subconverges in the sense described above.
\endproclaim
\proclabel{CompactnessPSLC}
\newsubhead{Degenerate submanifolds}[DegenerateSubmanifolds]
Recall now that, whereas the compactness result of Theorem \procref{CompactnessPSLC} concerns special legendrian submanifolds of $M$, what we actually require is a compactness result for constant scalar curvature, spacelike hypersurfaces in $\Bbb{R}^{d,1}$. Bearing in mind Lemma \procref{Dictionary}, such a compactness result will follow once we have identified all positive special legendrian submanifolds of $M$ which do not project down to constant scalar curvature, spacelike hypersurfaces in $\Bbb{R}^{d,1}$. However, these are precisely the positive special legendrian submanifolds over which the restriction of $m$ is degenerate at some point. We now study how this property affects the global structure of such submanifolds.
\par
Consider first a point $p\in M$, and a lagrangian subspace $E$ of the fibre $\alpha_p$. We say that $X\in E$ is an {\sl eigenvector} of $m$ over $E$ with eigenvalue $\lambda$ whenever
$$
m(X,Y) = \lambda g(X,Y)
$$
for all other $Y$ in $E$. In particular, when $0$ is an eigenvalue of $m$ over $E$, we define the {\sl nullity} of $m$ over $E$ to be the multiplicity of this eigenvalue, and we define it to be equal to $0$ otherwise. Consider now a positive, special legendrian submanifold $\hat{\Sigma}$ of $M$ so that, in particular, $T_p\hat{\Sigma}$ is a lagrangian subspace of the fibre $\alpha_p$ at every point $p\in\hat{\Sigma}$. In this and the following section, we will show that, in the case of interest to us, the nullity of $m$ is constant over $\hat{\Sigma}$. To this end, for all $p\in M$, for any lagrangian subspace $E$ of $\alpha_p$, and for all $1\leq k\leq d$, define
$$
\phi^k(E) := \minf_{(X_1,...,X_k)} \sum_{a=1}^k m(X_a,X_a),
$$
where $(X_1,...,X_k)$ ranges over all orthonormal $k$-tuples of vectors in $E$. By abuse of notation, for all $1\leq k\leq d$, define also $\phi^k:\hat{\Sigma}\rightarrow[0,\infty[$ by
$$
\phi^k(p) := \phi^k(T_p\Sigma).
$$
The following lemma will prove useful.
\proclaim{Lemma \nextprocno}
\noindent Let $p$ be a point in $M$ and let $E$ be a lagrangian subspace of the fibre $\alpha_p$. There exists an orthonormal basis of $E$ which simultaneously diagonalises both $m(\cdot,\cdot)$ and $m(\cdot,J\cdot)$.
\endproclaim
\proclabel{OrthogonalToMinimizers}
\proof Let $A,B:E\rightarrow E$ be the unique linear maps such that, for all $X,Y\in E$,
$$\eqalign{
m(X,Y) &= g(X,AY),\ \text{and}\cr
m(X,\opJ Y) &= g(X,BY).\cr}
$$
It suffices to show that $A$ and $B$ commute. Observe first that these linear maps are given by
$$\eqalign{
A &= \pi\opJ\opR,\ \text{and}\cr
B &= \pi\opJ\opR\opJ,}
$$
where $\pi:\alpha_p\rightarrow E$ here denotes the orthogonal projection with respect to $g$. Consequently,
$$\eqalign{
[A,B] &= \pi\opJ\opR\pi\opJ\opR\opJ - \pi\opJ\opR\opJ\pi\opJ\opR\cr
&= \pi\opR\opJ\pi\opR\opJ^2 + \pi\opR\opJ^2\pi\opJ\opR\cr
&= -\pi\opR(\opJ\pi+\pi\opJ)\opR.\cr}
$$
However, since $E$ is lagrangian, its orthogonal complement in $\alpha_p$ with respect to $g$ is given by
$$
E^\perp = \opJ E,
$$
so that
$$
\opJ\pi + \pi\opJ = \opJ,
$$
and
$$
[A,B] = -\pi\opR\opJ\opR = \pi\opR^2\opJ = \pi\opJ = 0,
$$
as desired.\qed
\proclaim{Lemma \nextprocno}
\noindent If $\hat{\Sigma}$ is positive then, for all $1\leq k\leq d$, there exists a continuous function $f^k$ over $\hat{\Sigma}$ such that
$$
\Delta\phi^{k} + f^k\phi^{k} \leq -\opTr(m)\sum_{a=1}^{k}\langle X_{a,i},X_{a,i}\rangle m(X_a,\opJ X_a)\eqnum{\nexteqnno[LaplacianOfPhi]}
$$
in the viscosity sense, where $(X_1,...,X_{k})$ is any orthonormal $k$-tuple of vectors in $E$ realising $\phi^{k}$ and, for each $a$, $X_a=:(X_{a,r},X_{a,i})$.
\endproclaim
\proclabel{FormulaForDeltaPhi}
\proof Choose $p\in\hat{\Sigma}$. Bearing in mind Lemma \procref{OrthogonalToMinimizers}, let $e_1,...,e_d$ be an orthonormal basis of joint eigenvectors of $m(\cdot,\cdot)$ and $m(\cdot,\opJ\cdot)$ over $T_p\hat{\Sigma}$, chosen in such a manner that $m(e_1,e_1)\leq...\leq m(e_d,e_d)$. We extend $e_1,...,e_d$ to an orthonormal frame field over $\hat{\Sigma}$ in a neighbourhood of $p$ which is parallel along geodesics leaving $p$. In particular,
$$
\phi^{k}(p) = \sum_{a=1}^{k}m(e_a(p),e_a(p)),
$$
and for all other $q$ near $p$,
$$
\phi^{k}(q)\leq \sum_{a=1}^{k}m(e_a(q),e_a(q)).
$$
Choose $1\leq a\leq k$ and, for ease of presentation, set $X:=e_a$. Consider now the function $m(X,X)$. We have
$$\eqalign{
\Delta m(X,X)
&= \sum_{b=1}^d e_be_b m(X,X)\cr
&= \sum_{b=1}^d 2e_b m\left(\overline{\nabla}_{e_b}X,X\right)\cr
&= \sum_{b=1}^d 2m\left(\overline{\nabla}_{e_b}\overline{\nabla}_{e_b}X,X\right) + 2m\left(\overline{\nabla}_{e_b}X,\overline{\nabla}_{e_b}X\right)\cr
&= \sum_{b=1}^d 2m\left(\nabla_{e_b}\nabla_{e_b}X + \hat{A}_{e_b}\nabla_{e_b}X + \overline{\nabla}_{e_b}\hat{A}_{e_b}X,X\right) + 2m\left(\overline{\nabla}_{e_b}X,\overline{\nabla}_{e_b}X\right).\cr}
$$
Since $X$ is parallel along geodesics leaving $p$, this yields
$$
\Delta m(X,X) = \sum_{b=1}^d 2m\left(\overline{\nabla}_{e_b}\hat{A}_{e_b}X,X\right) + 2m\left(\hat{A}_{e_b}X,\hat{A}_{e_b}X\right).
$$
Since $\hat{\Sigma}$ is positive, $m$ is non-positive over $\opJ T_p\hat{\Sigma}$, and so
$$
\Delta m(X,X) \leq \sum_{b=1}^d 2m\left(\overline{\nabla}_{e_b}\hat{A}_{e_b}X,X\right).
$$
Since $m(\cdot,\cdot)$ and $m(\cdot,\opJ\cdot)$ are both diagonal with respect to the basis $e_1,...,e_d$, this becomes
$$
\Delta m(X,X) \leq \sum_{b=1}^d 2g\left(\overline{\nabla}_{e_b}\hat{A}_{e_b}X,X\right)m(X,X) + 2g\left(\overline{\nabla}_{e_b}\hat{A}_{e_b}X,\opJ X\right)m(\opJ X,X).
$$
Since $g\left(\hat{A}_YX,Z\right)$ and $g\left(\hat{A}_YX,\opJ \hat{A}_WZ\right)$ always vanish, this yields
$$
\Delta m(X,X)\leq\sum_{b=1}^d (-2)g\left(\hat{A}_{e_b}X,\hat{A}_{e_b}X\right) m(X,X) + 2(\nabla_{e_b}\opII)(e_b,X,X)m(\opJ X,X).
$$
Finally, applying the Codazzi-Mainardi equations \eqnref{CodazziMainardi} yields
$$\eqalign{\Delta m(X,X)&\leq\sum_{b=1}^d\bigg[(-2)g\left(\hat{A}_{e_b}X,\hat{A}_{e_b}X\right) m(X,X)\cr
&\qquad\qquad + 2(\nabla_X\opII)(e_b,e_b,X)m(\opJ X,X)\phantom{\bigg]}\cr
&\qquad\qquad\qquad\qquad + 2g(\overline{R}_{e_b X}e_b,\opJ X) m(\opJ X,X)\bigg].\cr}\eqnum{\nexteqnno[FormulaForDeltaXX]}
$$
Since $m$ is non-negative over $T_p\hat{\Sigma}$, the first term on the right hand side is absorbed into $f^k\phi^k$, and since the second vanishes by Lemma \procref{Minimal}, it only remains to study the contribution of the third. However, by \eqnref{CurvatureOfContactBundle},
$$\eqalign{
\sum_{b=1}^dg(\overline{R}_{e_b X}e_b,\opJ X) m(\opJ X,X)
&=\sum_{b=1}^d(-1)\bigg[\langle(\opJ X)_r,e_{b,i}\rangle\langle e_{b,r},X_i\rangle\cr
&\qquad+\langle(\opJ X)_i,e_{b,i}\rangle\langle e_{b,i},X_i\rangle\bigg]m(\opJ X,X)\cr
&\qquad\qquad+\bigg[\langle e_{b,r},e_{b,i}\rangle\langle(\opJ X)_r,X_i\rangle\cr
&\qquad\qquad\qquad+\langle e_{b,i},e_{b,i}\rangle\langle(\opJ X)_i,X_i\rangle\bigg]m(\opJ X,X)\cr}
$$
Thus, using again the fact that $e_1,...,e_d$ diagonalises both $m(\cdot,\cdot)$ and $m(\cdot,\opJ\cdot)$,
$$\eqalign{
\sum_{b=1}^dg(\overline{R}_{e_b X}e_b,\opJ X) m(\opJ X,X)
&=\bigg[\langle X_i,X_i\rangle\langle X_r,X_i\rangle\cr
&\qquad-\langle X_r,X_i\rangle\langle X_i,X_i\rangle\bigg]m(\opJ X,X)\cr
&\qquad\qquad-\sum_{b=1}^d\langle e_{b,r},e_{b,i}\rangle\langle X_i,X_i\rangle m(\opJ X,X)\cr
&\qquad\qquad\qquad+\sum_{b=1}^d\langle e_{b,i},e_{b,i}\rangle\langle X_r,X_i\rangle m(\opJ X,X),\cr}
$$
so that
$$\eqalign{
\sum_{b=1}^dg(\overline{R}_{e_b X}e_b,\opJ X) m(\opJ X,X)
&=-\frac{1}{2}\opTr(m)\langle X_i,X_i\rangle m(\opJ X,X)\cr
&\qquad+\frac{1}{2}\sum_{b=1}^d\langle e_{b,i},e_{b,i}\rangle m(\opJ X,X)m(X,X).\cr}
$$
The second term on the right hand side is again absorbed into $f^k\phi^k$, and the result now follows upon summing over all $a$.\qed
\newsubhead{The refined special lagrangian angle}[RefiningTheSpecialLegendrianCondition]
We now show how positive special legendrian submanifolds divide into disjoint families, one of which will have the properties that we require. Indeed, consider a point $p\in M$, let $\Cal{R}$ and $\Cal{I}$ be respectively the real and imaginary subspaces of the fibre $\alpha_p$ as defined in \eqnref{RealAndImaginaryComponents}, and let $E$ be a positive, lagrangian subspace of the fibre $\alpha_p$. Upon perturbing $E$ slightly if necessary, we may suppose that it is the graph over $\Cal{R}$ of some symmetric matrix $A:\Cal{R}\rightarrow\Cal{I}$. Bearing in mind \eqnref{GraphicalSLCondition}, we see that $E$ is $\theta$-special lagrangian whenever
$$
\Theta(E) := \sum_{i=1}^d\opArctan(\lambda_i) = -\theta\ \text{mod}\ \pi\Bbb{Z},
$$
where $0<\lambda_1\leq...\leq\lambda_d$ are the eigenvalues of $A$. The function $\Theta$ extends continuously to all positive special lagrangian subspaces, including those that are not graphs. We call $\Theta$ the {\sl refined special lagrangian angle} of $E$. Of particular interest to us will be the case where $\Theta=\pi/2$.
\proclaim{Lemma \nextprocno}
\noindent Let $E$ be a positive lagrangian subspace of the fibre $\alpha_p$ with refined special lagrangian angle equal to $\pi/2$. For all $1\leq k\leq (d-1)$, there exists an orthonormal $k$-tuple $(X_1,...,X_k)$ of unit vectors in $E$ realising $\phi^k(E)$ such that, for all $a$,
$$
m(X_a,\opJ X_a)\geq0.
$$
\endproclaim
\proclabel{MJIsPositive}
\proof Indeed, upon perturbing $E$ slightly if necessary, we may suppose that it is the graph of a positive definite, symmetric matrix $A$. If $0<\lambda_1\leq...\leq\lambda_d$ are the eigenvalues of this matrix, then
$$
0 < \sum_{a=1}^d\opArctan(\lambda_a) = \Theta(E) = \pi/2.
$$
In particular,
$$
\opArctan(\lambda_{d-1}) \leq \frac{\pi}{4},
$$
and
$$
\opArctan(\lambda_{d-1}) \leq \opArctan(\lambda_d) \leq \frac{\pi}{2} - \opArctan(\lambda_{d-1}),
$$
so that
$$
\lambda_{d-1} \leq 1,
$$
and
$$
\lambda_{d-1}\leq\lambda_d\leq\frac{1}{\lambda_{d-1}}.
$$
Furthermore, $\lambda_d=1/\lambda_{d-1}$ if and only if $d=2$ and $\lambda_1=\lambda_2=1$. Suppose now that $(X_1,...,X_k)$ is an orthonormal $k$-tuple of vectors in $E$ which realises $\phi^k(E)$. Then, for all $1\leq a\leq k$, we may suppose that
$$
X_a = \frac{1}{1+\mu_a^2}(X_a',\mu_a X_a'),
$$
for some unit eigenvector $X_a'$ of $A$ with eigenvalue $\mu_a$. However, for each $a$,
$$
m(X,X) = \frac{2\mu_a}{1+\mu_a^2}.
$$
On the other hand, a straightforward calculation shows that
$$
\frac{2\lambda_d}{1 + \lambda_d^2} \geq \frac{2\lambda_{d-1}}{1 + \lambda_{d-1}^2} = \frac{2/\lambda_{d-1}}{1 + 1/\lambda_{d-1}^2},
$$
and for all $1\leqslant b\leqslant (d-1)$,
$$
\frac{2\lambda_b}{1+\lambda_b^2} \leq \frac{2\lambda_{d-1}}{1+\lambda_{d-1}^2}.
$$
It follows that
$$
\mu_a\leq\lambda_{d-1}\leq 1,
$$
so that
$$\eqalign{
m(X,\opJ X) &= \frac{1}{1 + \mu_a^2}m((X_a',\mu_a X_a'),(-\mu_a X_a',X_a'))\cr
&= \frac{1-\mu_a^2}{1+\mu_a^2}\cr
&\geq0,}
$$
as desired.\qed
\proclaim{Theorem \nextprocno}
\noindent If $\hat{\Sigma}$ is a positive, $\pi/2$-special legendrian submanifold of $M$ with refined special lagrangian angle equal to $\pi/2$, then $m$ has constant nullity over $\hat{\Sigma}$.
\endproclaim
\proclabel{TotalDegeneracy}
\proof First choose $1\leq k\leq (d-1)$. By Lemmas \procref{FormulaForDeltaPhi} and \procref{MJIsPositive}, there exists a continuous function $f^k:\hat{\Sigma}\rightarrow\Bbb{R}$ such that
$$
\Delta\phi^{k} + f^k\phi^{k} \leq 0.
$$
It follows by the strong maximum principle (c.f. Theorem $3.5$ of \cite{GilbTrud}, and the subsequent discussion) that if $\phi^k$ vanishes at a single point, then it vanishes identically. We conclude that either $m$ has constant nullity over $\hat{\Sigma}$, or that its nullity is at least $(d-1)$ at every point of this submanifold. Finally, since $\hat{\Sigma}$ is special legendrian with $\Theta=\pi/2$, the nullity of $m$ cannot be equal to $(d-1)$, so that the latter case only occurs when $m$ vanishes identically. This completes the proof.\qed
\newsubhead{The geometry of curtain submanifolds}[TheGeometryOfCurtainSubmanifolds]
Following \cite{LabA}, we call {\sl curtain submanifolds} those positive, special legendrian submanifolds over which the restriction of $m$ has non-trivial nullity. We now describe their geometry. Observe first that, for all $d'<d$, the product $\Bbb{R}^{d-d'}\times M^{d'}$ naturally embeds into $M^d$ and, if $\hat{\Sigma}$ is a positive special legendrian submanifold of $M^{d'}$, then $\Bbb{R}^{(d-d')}\times\hat{\Sigma}$ is a positive special legendrian submanifold of $M^d$ with the same refined special lagrangian angle. We now show that, up to the action of an element of $\opO(d,1)$, such products account for all curtain submanifolds. To this end, the following construction will prove useful. Let $\Bbb{H}^d$ denote the upper component of the unit pseudosphere in $\Bbb{R}^{d,1}$, that is,
$$
\Bbb{H}^d = \left\{ y\ |\ \|y\|^2 = -1,\ y_{d+1}>0\right\}.
$$
Observe that $\Bbb{H}^d$ is isometric to $d$-dimensional hyperbolic space. Consider now an embedded submanifold $S$ in $\Bbb{H}^d$ and denote its normal bundle in $\Bbb{H}^d$ by $NS$, that is
$$
NS := \left\{ (x,y)\ |\ y\in S,\ x\perp y,\ x\perp T_yS\right\}.
$$
Given a smooth function $\xi:S\rightarrow\Bbb{R}^{d,1}$, we now define
$$
N^\xi S := \left\{ (x+\xi(y),y)\ |\ y\in S, x\in N_yS\right\}.
$$
In particular, we may assume that $\xi$ is everywhere normal to $N_yS$. We now establish under which conditions $N^\xi S$ is positive special legendrian.
\proclaim{Lemma \nextprocno}
\noindent Let $\xi:S\rightarrow\Bbb{R}^{d,1}$ be such that $\xi(y)$ is orthogonal to $N_yS$ for all $y\in S$.
\medskip
\myitem{(1)} $N^\xi S$ is a legendrian submanifold of $M^d$ if and only if
$$
\xi(y) = \nabla\phi(y) - \phi(y)y
$$
for some smooth function $\phi:S\rightarrow\Bbb{R}$;
\medskip
\myitem{(2)} $N^\xi S$ is positive if and only if $S$ is totally geodesic; and
\medskip
\myitem{(3)} $N^\xi S$ is positive special legendrian if and only if
$$
\opIm(e^{i\theta}\opDet(\opHess(\phi) + (i-\phi)\opId)) = 0,
$$
for some $\theta\in\Bbb{R}$, where $\phi:S\rightarrow\Bbb{R}$ is the function defined in $(1)$.
\endproclaim
\proclabel{DescribingCurtainSubmanifolds}
\remark In particular, if $N^\xi S$ is positive special legendrian then, upon applying an element of $\opO(d,1)$, we may suppose that $S$ is an open subset of the totally geodesic subspace
$$
(\left\{0\right\}\times\Bbb{R}^{d',1})\minter\Bbb{H}^d.
$$
It then follows that
$$
N^\xi S = \Bbb{R}^{(d-d')}\times\hat{\Sigma}
$$
for some $d'<d$ and for some $d'$-dimensional submanifold $\hat{\Sigma}$ of $M^{d'}$. In particular, it is straightforward to verify that $\hat{\Sigma}$ is also positive special legendrian, so that $N^\xi S$ is one of the products described above.
\medskip
\proof Let $0\leq m<d$ be the dimension of $S$. Consider a point $(x,y)\in NS$. Let $e_1,...,e_m$ be an orthonormal basis of $T_yS$ and extend this to an orthonormal basis $e_1,...,e_d$ of $T_y\Bbb{H}^d$. In particular $e_{m+1},...,e_d$ is an orthonormal basis of $N_y S$. A basis of $T_{(x+\xi(y),y)}N^\xi S$ is now given by
$$
f_i := \left\{
\matrix (A(x,y)e_i + D_{e_i}\xi(y),e_i)\hfill&\ \text{if}\ 1\leq i\leq m,\ \text{and}\hfill\cr
(e_i,0)\hfill&\ \text{if}\ m<i\leq d,\hfill\cr\endmatrix\right.
$$
where, for all $(x,y)\in NS$, $A(x,y):T_yS\rightarrow T_yS$ here denotes the shape operator of $S$ at the point $y$ with respect to the normal vector $x$. For all $1\leq i,j\leq m$, denote
$$
a_{ij}(x,y):=\langle A(x,y)e_i + D_{e_i}\xi(y),e_j\rangle.
$$
Suppose now that $N^\xi S$ is legendrian. In particular, $T_{(x+\xi(y),y)}N^\xi S$ is contained in the fibre $\alpha_{(x+\xi(y),y)}$ so that, for all $i$,
$$
\langle D_{e_i}\xi(y),y\rangle = 0.
$$
Thus, setting $\phi:=\langle\xi(y),y\rangle$, we have
$$
\langle\nabla\phi(y),e_i\rangle=D_{e_i}\langle \xi(y),y\rangle = \langle \xi(y),D_{e_i}y\rangle = \langle\xi(y),e_i\rangle,
$$
and since $y$ has negative norm-squared, it follows that
$$
\xi(y) = \nabla\phi(y) - \phi(y)y.
$$
Conversely, if $\xi$ has the above form, then we readily verify that $T_{(x+\xi(y),y)}N^\xi S$ is contained in the fibre $\alpha_{(x+\xi(y),y)}$ and that $a_{ij}$ is symmetric, so that $N^\xi S$ is legendrian. This proves $(1)$. Observe now that $N^\xi S$ is positive if and only if the matrix $a_{ij}$ is positive definite for all $x$ and for all $y$. However, since $A$ is linear in $x$, this holds if and only if $A(x,y)$ vanishes identically, and this proves $(2)$. Finally, observe that the hessian of $\phi$ over $\Sigma$ is given by
$$
\opHess(\phi)(y)(e_i,e_j) = \langle D_{e_i}\nabla\phi(y),e_j\rangle = \langle D_{e_i}\xi(y),e_j\rangle + \phi\langle e_i,e_j\rangle.
$$
However, as in \eqnref{GraphicalSLCondition}, $N^\xi S$ is special legendrian if and only if
$$
\opIm(e^{i\theta}\opDet(a + i\opId)) = 0,
$$
for some $\theta$, and since $A(x,y)$ vanishes identically, this proves $(3)$.\qed
\proclaim{Theorem \nextprocno}
\noindent If $\hat{\Sigma}$ is a complete curtain submanifold of $M^d$ with refined special lagrangian angle equal to $\pi/2$ then, up to the action of an element of $\opO(d,1)$,
$$
\hat{\Sigma} = \Bbb{R}^{(d-d')}\times\hat{\Sigma}',
$$
where $d'<d$ and $\hat{\Sigma}'$ is a complete, positive, special legendrian submanifold of $M^{d'}$ with refined special lagrangian angle equal to $\pi/2$.
\endproclaim
\proclabel{CurtainSubmanifolds}
\remark Observe that the only positive special legendrian submanifolds of $M^1$ with refined special lagrangian angle equal to $\pi/2$ are simply the fibres of this bundle over $\Bbb{R}^{1,1}$. In this manner, we recover the structure of $2$-dimensional curtain submanifolds studied by Labourie in \cite{LabA}.
\medskip
\proof By Lemma \procref{DescribingCurtainSubmanifolds}, it suffices to show that every point $p$ of $\hat{\Sigma}$ has a neighbourhood $\Omega$ of the form
$$
\Omega = N^\xi S,
$$
for some embedded submanifold $S$ of $\Bbb{H}^d$ and some smooth function $\xi:S\rightarrow\Bbb{R}^{d,1}$. To this end, observe first that
$$
\opKer(m) = (T\hat{\Sigma}\minter\Cal{R})\oplus(T\hat{\Sigma}\minter\Cal{I}).
$$
Furthermore, since $\hat{\Sigma}$ has refined special lagrangian angle equal to $\pi/2$, the distribution $T\hat{\Sigma}\minter\Cal{I}$ has dimension at most $1$, and is non-trivial if and only if $m$ vanishes identically over $\hat{\Sigma}$. In particular, the distribution $T\hat{\Sigma}\minter\Cal{R}$ has constant non-zero dimension. Observe now that $T\hat{\Sigma}\minter\Cal{R}$ also coincides with $\opKer(D\pi)$, where $\pi:M^d\rightarrow\Bbb{H}^d$ is the projection onto the second factor. However, considered as a $1$-form taking values in $\Bbb{R}^{d,1}$, we have
$$
dD\pi=0,
$$
and it follows that $T\hat{\Sigma}\minter\Cal{R}$ is integrable. Let $\Cal{F}$ denote the smooth foliation of $\hat{\Sigma}$ that it defines, and observe that $\pi$ is constant over every leaf of this foliation.
\par
Consider now a point $p\in\hat{\Sigma}$, let $\hat{S}$ be a smooth submanifold passing through $p$ which is transverse to $\Cal{F}$, and let $S$ be its image under the projection $\pi$. Upon reducing $\hat{S}$ if necessary, we may suppose that $S$ is an embedded submanifold of $\Bbb{H}^d$, and that $\pi$ restricts to a diffeomorphism of $\hat{S}$ onto $S$. Let $\Omega$ be the neighbourhood of $p$ consisting of the union of all leaves of $\Cal{F}$ which pass through $\hat{S}$. We show that $\Omega$ has the desired form. Indeed, let $\hat{q}$ be another point of $\Omega$, denote $q:=\pi(\hat{q})$, let $L$ be the leaf of $\Cal{F}$ passing through $\hat{q}$, and let $X:=(X_r,0)$ be a tangent vector to $L$ at this point. Since $T_{\hat{q}}\hat{\Sigma}$ is contained in $\alpha_{\hat{q}}$, we have
$$
\langle X_r,q\rangle = 0,
$$
and since $m$ is non-negative semi-definite over $T_{\hat{q}}\hat{\Sigma}$, for all other tangent vectors $Y:=(Y_r,Y_i)$ to $\hat{\Sigma}$ at $\hat{q}$, we have
$$
m(X,Y)^2 \leq m(X,X)m(Y,Y) = 0,
$$
so that, for all tangent vectors $Y_i$ to $S$ at $q$, we have
$$
\langle X_r,Y_i\rangle = 0.
$$
It follows that
$$
T_{\hat{q}} = N_qS,
$$
and since $\hat{\Sigma}$ is complete, upon integrating we find that that there exists a point $\xi\in\Bbb{R}^{d,1}$ such that
$$
L_{\hat{q}} = \xi + N_qS,
$$
as desired.\qed
\newhead{Hypersurfaces in Minkowski spacetimes}[HypersurfacesInMinkowskiSpace]
\newsubhead{Affine deformations and GHMC Minkowski spacetimes}[QuasifuchsianRepresentations]
Let $\opPSO(d,1)$ denote the group of linear isometries of $\Bbb{R}^{d,1}$ which preserve the spatial and temporal orientations. Consider now the upper unit pseudo-sphere in $\Bbb{R}^{d,1}$,
$$
\Bbb{H}^d := \left\{ x\in\Bbb{R}^{d,1}\ |\ \|x\|^2 = -1,\ x_{d+1} > 0\right\},
$$
and recall that the semi-riemannian metric of $\Bbb{R}^{d,1}$ restricts to a complete hyperbolic metric over this submanifold. Since $\opPSO(2,1)$ preserves $\Bbb{H}^d$, it identifies in this manner with the group of orientation preserving isometries of hyperbolic space. A subgroup $\Gamma$ of $\opPSO(d,1)$ is said to be {\sl kleinian} whenever it is discrete, cocompact and torsion free (c.f. \cite{KapovichII}). In particular, the quotient of $\Bbb{H}^d$ by a kleinian subgroup is a compact hyperbolic manifold and, conversely, the fundamental group of any compact hyperbolic manifold identifies with some kleinian subgroup.
\par
The group of affine isometries of $\Bbb{R}^{d,1}$ which preserve the spatial and temporal orientations is given by the semidirect product $\opPSO(d,1)\ltimes\Bbb{R}^{d,1}$, where the group law is given by
$$
(\alpha,x)\cdot(\beta,y) := (\alpha\beta,x + \alpha(y)).
$$
Given a kleinian subgroup $\Gamma$, an {\sl affine deformation} is a homomorphism $\rho:\Gamma\rightarrow\opPSO(d,1)\ltimes\Bbb{R}^{d,1}$ of the form
$$
\rho(\alpha)=(\alpha,\tau(\alpha)),
$$
for some map $\tau:\Gamma\rightarrow\Bbb{R}^{d,1}$. By abuse of notation, we denote the image of $\rho$ in $\opPSO(d,1)\ltimes\Bbb{R}^{d,1}$ by $\Gamma\ltimes\tau$. The homomorphism property of $\rho$ is equivalent to the {\sl cocycle condition} on $\tau$, namely
$$
\tau(\alpha\beta) = \tau(\alpha) + \alpha\tau(\beta),
$$
for all $\alpha,\beta\in\Gamma$. It follows that the set of affine deformations of $\Gamma$ naturally identifies with the set of $\Bbb{R}^{d,1}$-cocycles over this group, which itself constitutes a finite-dimensional vector space.
\par
GHMC Minkowski spacetimes are parametrised by affine deformations as follows. First, recall that the future cone in $\Bbb{R}^{d,1}$ is given by
$$
C := \left\{ x\in\Bbb{R}^{d,1}\ |\ \|x\|^2 \leq 0,\ x_{d+1} \geq 0\right\}.
$$
We say that a closed, convex subset $K$ of $\Bbb{R}^{d,1}$ is future-complete whenever
$$
K+C:=\left\{x+y\ |\ x\in K,\ y\in C\right\}\subseteq K.
$$
In \cite{Bonsante}, building on the work \cite{Mess} of Mess, Bonsante shows that, for any affine deformation $\Gamma\ltimes\tau$ in $\opPSO(d,1)\ltimes\Bbb{R}^{d,1}$, there exists a unique future-complete, closed, convex subset $K$ of $\Bbb{R}^{d,1}$ such that
\medskip
\myitem{(1)} $K$ is invariant under the action of $\Gamma\ltimes\tau$,
\medskip
\myitem{(2)} $\Gamma\ltimes\tau$ acts properly discontinuously on the interior of $K$, and
\medskip
\myitem{(3)} $K$ is maximal with respect to inclusion amongst all future-complete, closed, convex subsets of $\Bbb{R}^{d,1}$ which satisfy $(1)$ and $(2)$.
\medskip
\noindent The quotient $K^o/\Gamma\ltimes\tau$ is a GHMC Minkowski spacetime. Throughout the sequel, we will refer to the set $K$ as {\sl Bonsante's domain} for $\tau$, and we will refer to GHMC spacetimes that arise in this manner as {\sl Bonsante spacetimes}. In \cite{Barbot}, Barbot shows that, up to reversal of the temporal orientation, every GHMC Minkowski spacetime which is not a translation spacetime or a Misner spacetime is, up to a finite cover, a twisted product of a Bonsante spacetime with a euclidean torus. In particular, the general result readily follows from the result for Bonsante spacetimes, and we therefore restrict our attention henceforth to this case.
\par
Bonsante's construction can also be interpreted as follows. Given a fixed kleinian subgroup $\Gamma$, there exists a set valued function $K$ which maps every $\Bbb{R}^{d,1}$-cocycle $\tau$ over $\Gamma$ to Bonsante's domain $K(\tau)$ for $\tau$. It can be deduced from \cite{Bonsante} that this function is continuous with respect to the local Hausdorff topology, that is, if $(\tau_m)$ is a sequence of cocycles converging to $\tau_\infty$ then, for all $r>0$,
$$
K(\tau_m)\minter B_r(0)\rightarrow K(\tau_\infty)\minter B_r(0)
$$
in the Hausdorff sense, where $B_r(0)$ here denotes the open ball of (euclidean) radius $r$ about the origin in $\Bbb{R}^{d,1}$. Furthermore, all supporting hyperplanes to $K(\tau)$ are either spacelike or null and the intersection of $K(\tau)$ with any spacelike hyperplane is compact. It follows, in particular, that if $(x_m)$ is a divergent sequence of boundary points of $K(\tau)$ and if, for all $m$, $N_m$ is a unit, timelike supporting normal at the point $x_m$, then the sequence $(N_m)$ also approaches the light cone as $m$ tends to infinity.
\par
At this stage, it is worthwhile to verify the existence of non-trivial examples of $\Bbb{R}^{d,1}$-cocycles and affine deformations. Indeed, although, in the $(2+1)$-dimensional case, large families of cocycles are obtained via the natural identification of $\Bbb{R}^{2,1}$-cocyles with tangent vectors to the Teichm\"uller space of the surface $\Bbb{H}^2/\Gamma$ (c.f. \cite{Goldman}), in higher dimensions, Cartan-Weyl local rigidity makes the problem of constructing non-trivial examples a good deal more subtle. One nice technique, however, involves interpreting cocycles as infinitesimal variations of the hyperbolic manifold $\Bbb{H}^d/\Gamma$ within the space of flat conformal manifolds. Indeed, consider the canonical injection
$$
\theta_0:\Gamma\rightarrow\opPSO(d,1)\rightarrow\opPSO(d+1,1),
$$
and recall that any homomorphism $\theta:\Gamma\rightarrow\opPSO(d+1,1)$ which is sufficiently close to $\theta_0$ is the holonomy of some flat conformal structure over the manifold $\Bbb{H}^d/\Gamma$ (c.f. \cite{KulkPink} and \cite{SmiFCS}). Suppose now that $\theta_0$ extends to a non-trivial smooth family $(\theta_t)_{t\in]-\epsilon,\epsilon[}$ of homomorphisms of $\Gamma$ into $\opPSO(d+1,1)$. Then $\tilde{\tau}:=\theta_0^{-1}(\partial_t\theta)_0$ defines a cocycle taking values in the Lie algebra $\frak{so}(d+1,1)$. However, when the space $\frak{so}(d+1,1)$ is considered as a representation of $\opPSO(d,1)$ via the adjoint action, it naturally decomposes as
$$
\frak{so}(d+1,1) = \frak{so}(d,1)\oplus\Bbb{R}^{d,1},
$$
and since the first component of $\tilde{\tau}$ vanishes by Cartan-Weyl local rigidity, this map is entirely determined by its second component, which we readily verify to be the desired cocycle in $\Bbb{R}^{d,1}$. Finally, there are various known constructions of non-trivial, smooth families of flat conformal structures. The simplest involves bending a compact hyperbolic manifold around a totally geodesic hypersurface, whenever such a hypersurface exists (c.f. \cite{JohnsonMillson} and \cite{Kourouniotis}). We refer the reader to \cite{Apanasov}, \cite{Kapovich} and \cite{Scannell} for more details of this and other constructions.
\par
We conclude this section by studying the case where the cocycle vanishes, which is known as the {\sl fuchsian} case, and forms the starting point of our construction. In particular, Theorem \procref{MainTheorem} trivially holds in this case. Indeed, observe first that Bonsante's domain for the trivial cocycle $\tau=0$ coincides with the future cone $C$. Now, given $k\in]0,\infty[$, consider the hypersurface
$$
\Sigma_k := \left\{ x\in C\ |\ \|x\|^2 = -1/k^2\right\},
$$
and observe that this hypersurface has constant scalar curvature equal to $(-k^2)$. Furthermore, the family $(\Sigma_k)_{k\in]0,\infty[}$ constitutes a smooth foliation of the interior of $C$. Indeed, if we define the smooth submersion $\phi:C^o\rightarrow]0,\infty[$ by
$$
\phi(y) := \frac{-1}{\|y\|^2},
$$
then, for all $k$, $\Sigma_k$ is the level subset of $\phi$ at height $k$. Observe finally that $(\Sigma_k)$ tends to $\partial C$ in the local Hausdorff sense as $k$ tends to $+\infty$. This foliation, which we henceforth refer to as the {\sl fuchsian foliation}, will be of use at various stages in the sequel.
\newsubhead{Stability}[Stability]
Let $\Gamma\subseteq\opPSO(3,1)$ be a kleinian subgroup, let $\tau:\Gamma_0\rightarrow\Bbb{R}^{3,1}$ be a cocycle, and let $K:=K(\tau)$ be Bonsante's domain for $\tau$. For $k>0$, let $\Sigma_k\subseteq K$ be a spacelike, locally strictly convex hypersurface in $K$ of constant scalar curvature equal to $(-k^2)$ and invariant under the action of $\Gamma\ltimes\tau$. In this section, we study infinitesimal perturbations of $\Sigma_k$ corresponding to infinitesimal variations of the cocycle $\tau$.
\par
First recall that, by the Gauss-Codazzi equations, the scalar curvature of $\Sigma_k$ is given by
$$
S:=-\frac{1}{3}(\lambda_1\lambda_2 + \lambda_1\lambda_3+\lambda_2\lambda_3),
$$
where $\lambda_1$, $\lambda_2$ and $\lambda_3$ are its principal curvatures. The Jacobi operator of $\Sigma_k$ then measures the infinitesimal variation of scalar curvature resulting from an infinitesimal normal perturbation of this surface. More formally, let $N:\Sigma_k\rightarrow\Bbb{R}^{3,1}$ be the future-oriented, unit normal vector field over $\Sigma_k$, and for $f\in C_0^\infty(\Sigma_k)$, and for $t\in\Bbb{R}$, define
$$
\Phi_{f,t}(x) := x + tf(x)N(x).
$$
For any given $x$, and for sufficiently small $t$, $\Phi_{f,t}$ is also an embedding near $x$. In particular, letting $S_{f,t}(x)$ denote its scalar curvature at the point $x$, we define
$$
(Jf)(x) := \frac{\partial}{\partial t} S_{f,t}(x)\bigg|_{t=0},
$$
and we call $J$ the {\sl Jacobi operator} of $\Sigma_k$.
\par
Now let $A$ be the shape operator of $\Sigma_k$. Since $\Sigma_k$ is locally strictly convex, $A$ is everywhere positive definite. Define $B$ by
$$
B := \frac{1}{3}\left(\opTr(A)\opId - A\right).
$$
It is straightforward to see that $B$ is also everywhere positive definite, and that
$$
\opTr(BA) = 2k^2.
$$
\proclaim{Lemma \nextprocno}
\noindent The Jacobi operator of $\Sigma_k$ is given by
$$
Jf = B^{ij}A_{ij}^2f - B^{ij}f_{;ij},
$$
where the summation convention is implied and $f_{;ij}$ here denotes the hessian of $f$ along $\Sigma_k$.
\endproclaim
\proclabel{JacobiOperator}
\noindent Applying the maximum principle and the Fredholm alternative immediately yields
\proclaim{Corollary \nextprocno}
\noindent The Jacobi operator of $\Sigma_k$ is invertible.
\endproclaim
{\bf\noindent Proof of Lemma \procref{JacobiOperator}:\ }The following calculation is standard in the riemannian setting. We repeat it here in the lorentzian setting as care is required with signs that are different in certain places. Fix $f\in C^\infty(\Sigma_k)$. Define $I:\Sigma_k\times\Bbb{R}\rightarrow\Bbb{R}^{3,1}$ by
$$
I(x,t) = x + tf(x)N(x).
$$
Denote by $g$ the pull-back through $I$ of the Minkowski metric over $\Bbb{R}^{3,1}$. We extend $N$ to a vector field over $\Sigma_k\times]-\epsilon,\epsilon[$ such that, for all $(x,t)$, $N(x,t)$ is normal to the hypersurface, $\Sigma_k\times\left\{t\right\}$. Observe that, by definition, for all $x$
$$
f(x,0)N(x,0) = \partial_t.
$$
We claim that, along $\Sigma_k\times\left\{0\right\}$,
$$
\nabla_{\partial_t}N = \nabla^\Sigma f.
$$
Indeed,
$$\triplealign{
&\partial_t\langle N,N\rangle &= 0\cr
\Rightarrow&\langle\nabla_{\partial_t}N,N\rangle &= 0.\cr}
$$
Likewise, for any vector field $U$ tangent to $\Sigma_k$ and independent of $t$,
$$\triplealign{
&\partial_t\langle N,U\rangle &= 0\cr
\Rightarrow&\langle\nabla_{\partial_t}N,U\rangle &= -\langle N,\nabla_{\partial_t}U\rangle\cr
& &=-\langle N,\nabla_U\partial_t\rangle\cr
& &=-\langle N,\nabla_U fN\rangle\cr
& &=\langle\nabla^\Sigma f,U\rangle,\cr}
$$
where the last equality follows from the fact that $\langle N,N\rangle=-1$. This proves the assertion.
\par
We now define the endomorphism field $A$ over $\Sigma_k\times]-\epsilon,\epsilon[$ such that, for all $(x,t)$, $(AN)(x,t)=0$, and the restriction of $A(x,t)$ to the tangent space of $\Sigma_k\times\left\{t\right\}$ coincides with the shape operator of this hypersurface at this point. We are interested in the covariant derivative of $A$ in the time direction. Bearing in mind that $\Bbb{R}^{3,1}$ is flat, for any vector field $U$ tangent to $\Sigma_k$ and independent of time, we have
$$\eqalign{
(\nabla_{\partial_t}A)U
&=\nabla_{\partial_t}(AU) - A\nabla_{\partial_t}U\cr
&=\nabla_{\partial_t}\nabla_U N - A\nabla_U{\partial_t}\cr
&=\nabla_U\nabla_{\partial_t} N - A\nabla_U fN\cr
&=\nabla^\Sigma_U\nabla^\Sigma f - fA^2U.\cr}
$$
\par
Finally, the scalar curvature of $\Sigma_k$ is given by
$$
S = \sigma(A) := -\frac{1}{6}\left(\opTr(A)^2 - \opTr(A^2)\right).
$$
The derivative of the function $\sigma$ at $A$ is given by
$$
DS(A)M = -B^{ij}M_{ij},
$$
where the summation convention is implied, and the result now follows by the chain rule.\qed
\proclaim{Lemma \nextprocno}
\noindent Let $(\tau_t)_{t\in]-\epsilon,\epsilon[}$ be a smoothly varying family of $\Bbb{R}^{3,1}$-cocycles such that $\tau_0=\tau$. Upon reducing $\epsilon$ if necessary, there exists a unique, smoothly varying family $(\Sigma_{k,t})_{t\in]-\epsilon,\epsilon[}$ of spacelike, locally strictly convex hypersurfaces such that, for all $t$, $\Sigma_{k,t}$ is of constant scalar curvature equal to $(-k^2)$ and is invariant under the action of $\Gamma\ltimes\tau_t$.
\endproclaim
\proclabel{Openness}
\proof We first define a smooth family $(\Sigma_{k,t}')_{t\in]-\epsilon,\epsilon[}$ of spacelike, locally strictly convex hypersurfaces such that, for all $t$, the hypersurface $\Sigma_{k,t}'$ is invariant under the action of $\Gamma\ltimes\tau_t$ but does not necessarily satisfy the curvature condition. The desired family will then be obtained via a perturbation argument. In what follows, we identify $\Gamma$ with the action of $\Gamma\ltimes\tau_0$ over $\Sigma_0$. Let $N:\Sigma_k\rightarrow\Bbb{R}^{3,1}$ be the future-oriented, unit, normal vector field over $\Sigma_k$. Let $\phi\in C_0^\infty(\Sigma_k)$ be a smooth, positive function of compact support such that, for all $x\in\Sigma_k$,
$$
\sum_{\gamma\in\Gamma}\phi(\gamma(x)) = 1.
$$
Define $e:]-\epsilon,\epsilon[\times\Sigma_k\rightarrow\Bbb{R}^{3,1}$ by
$$
e_t(x) := e(x) + \sum_{\gamma\in\Gamma}\phi(\gamma(x))\left(\tau_t(\gamma^{-1})-\tau_0(\gamma^{-1})\right),
$$
and, for all $t$, denote $e_t:=e(t,\cdot)$. By construction, for all $t$, $e_t(\Sigma_k)$ is invariant under the action of $\Gamma\ltimes\tau_t$ and, upon reducing $\epsilon$ if necessary, we may suppose furthermore that it is embedded. This yields the desired family.
\par
For all $t$, let $N_t:\Sigma_k\rightarrow\Bbb{R}^{3,1}$ be the future-oriented, unit, normal vector field over $e_t$. Define $\Phi:C^\infty(\Sigma_k)\times]-\epsilon,\epsilon[\rightarrow C^\infty(\Sigma_k,\Bbb{R}^{3,1})$ by
$$
\Phi_{f,t}(x) = e_t(x) + f(x)N_t(x).
$$
Observe that if $f$ is $\Gamma$-invariant, then $\Phi_{f,t}$ is $\Gamma\ltimes\tau_t$ equivariant for all $t$. Furthermore, for sufficiently small $(f,t)$, $\Phi_{f,t}$ is an embedding. Let $C^\infty_\opinv(\Sigma_k)$ denote the space of smooth functions over $\Sigma_k$ that are $\Gamma$-invariant and define $S:C^\infty_\opinv(\Sigma_k)\times]-\epsilon,\epsilon[\rightarrow C^\infty_\opinv(\Sigma_k)$ such that, for all sufficiently small $(f,t)$, and for all $x$, $S_{f,t}(x)$ is the scalar curvature of the embedding $\Phi_{f,t}$ at the point $x$. Now, given $(l,\alpha)$, by defining the H\"older space $C^{l,\alpha}_\opinv(\Sigma_k)$ in a similar manner, we see that the functional $S$ extends continuously to a smooth map from $C^{l+2,\alpha}_\opinv(\Sigma_k)\times]-\epsilon,\epsilon[$ into $C^{l,\alpha}_\opinv(\Sigma_k)$. Furthermore, its partial derivative with respect to the first component at the point $(0,0)$ is simply the Jacobi operator $J$. Since $J$ is invertible, it follows by the implicit function theorem for maps between Banach manifolds that, upon reducing $\epsilon$ if necessary, there exists a smooth map $\phi:]-\epsilon,\epsilon[\rightarrow C^{l+2,\alpha}_\opinv(\Sigma_k)$ such that, for all $t$, $\Phi_{\phi(t),t}(\Sigma_k)$ has constant scalar curvature equal to $(-k^2)$. Furthermore, by elliptic regularity, these surfaces are smooth for all $t$, and this completes the proof.\qed
\newsubhead{Compactness}[Compactness]
\noindent We first obtain an elementary result concerning the position of a constant scalar curvature hypersurface inside a given GHMC Minkowski spacetime. As before, let $\Gamma\subseteq\opPSO(3,1)$ be a kleinian subgroup, let $\tau:\Gamma\rightarrow\Bbb{R}^{3,1}$ be a cocycle, and let $K:=K(\tau)$ be Bonsante's domain for $\tau$. For $r>0$, let $K^r$ denote the set of all points of $K$ lying at a (timelike) distance of at least $1/r$ from the boundary $\partial K$. For $k>0$, let $\Sigma_k\subseteq K$ be a spacelike, locally strictly convex hypersurface of constant scalar curvature equal to $(-k^2)$ which is invariant under the action of $\Gamma\ltimes\tau$, and let $\Sigma_k^+$ denote the future-complete, convex set bounded by $\Sigma_k$.
\proclaim{Lemma \nextprocno}
\noindent $K^k\subseteq\Sigma_k^+\subseteq K$.
\endproclaim
\proclabel{Position}
\proof Let $x$ be a point of $\partial K$. For all $\epsilon>0$, denote $x(\epsilon):=x+\epsilon e_4$, where $e_4$ here denotes the fourth canonical basis vector of $\Bbb{R}^{3,1}$, and let $C_{x(\epsilon)}$ be the future cone based on $x(\epsilon)$. Define $\phi:C_{x(\epsilon)}^o\rightarrow]0,\infty[$ by
$$
\phi(y) = \frac{-1}{\|y-x(\epsilon)\|^2},
$$
and observe that the level sets of $\phi$ are simply the leaves of the fuchsian foliation of $C_{x(\epsilon)}^o$ which was introduced in Section \subheadref{QuasifuchsianRepresentations}. In particular, for all $y\in C_{x(\epsilon)}^o$, the level set of $\phi$ passing through the point $y$ has constant scalar curvature equal to $\phi(y)$. Consider now a fundamental domain $\Omega$ of $\Sigma_k$ and observe that, for all $\epsilon>0$, there are only finitely many elements $\alpha$ of $\Gamma\ltimes\tau$ such that $\alpha(\Omega)$ has non-trivial intersection with $C_{x(\epsilon)}$. From this it follows that $\Sigma_k\minter C_{x(\epsilon)}$ is a relatively compact open subset of $\Sigma_k$, and so $\phi$ attains a minimum value at some point $y$, say, of this intersection. At this point, $\Sigma_k$ is an interior tangent to the leaf of curvature $\phi(y)$ so that, by the geometric maximum principle, $k^2\leq\phi(y)$. In particular, $\Sigma_k$ has trivial intersection with the set $\phi^{-1}(]0,k[)$, and the result follows by letting $\epsilon$ tend to $0$ and taking the union over all $x\in\partial K$.\qed
\medskip
Consider now a sequence $(k_m)$ of positive real numbers, and a sequence $(\tau_m)$ of cocycles. For all $m$, let $K_m:=K(\tau_m)$ be Bonsante's domain for $\tau_m$, and let $\Sigma_{k,m}\subseteq K_m$ be a spacelike, locally strictly convex hypersurface of constant scalar curvature equal to $(-k_m^2)$ which is invariant under the action of $\Gamma\ltimes\tau_m$. Suppose that $(k_m)$ and $(\tau_m)$ converge to $k_\infty$ and $\tau_\infty$ respectively. In particular $(K_m)$ converges in the local Hausdorff sense to $K_\infty:=K(\tau_\infty)$.
\proclaim{Lemma \nextprocno}
\noindent For every compact subset $X$ of $\Bbb{R}^{3,1}$, there exists a compact subset $Y$ of $\Bbb{H}^3$ such that, for any $m$, and for any point $x$ of $\Sigma_m\minter X$, if $N(x)$ is the future-oriented, unit normal vector of $\Sigma_m$ at this point, then $N(x)$ is an element of $Y$.
\endproclaim
\proclabel{UniformlySpacelike}
\proof Suppose the contrary. There exists a sequence $(x_m)$ such that, for all $m$, $x_m\in\Sigma_m\minter X$ but such that $(N_m(x_m))$ diverges. Since $\Gamma$ is cocompact, upon composing with suitable elements of $\Gamma\ltimes\tau_m$, we may suppose instead that the sequence $(N_m(x_m))$ remains within some fixed compact set, but that the sequence $(x_m)$ diverges. For all $m$, denote $\overline{K}_m:=K_m-x_m$, $\overline{K}_m^{k_m}:=K^{k_m}_m-x_m$ and $\overline{\Sigma}_m^+:=\Sigma_m^+-x_m$. Since $(K_m)$ converges towards $(K_\infty)$, we may suppose that both $(\overline{K}_m)$ and $(\overline{K}_m^{k_m})$ converge in the local Hausdorff sense to the future side of the same null hyperplane $H$, say. Furthermore, since $\overline{K}^{k_m}_m\subseteq\overline{\Sigma}_m^+\subseteq\overline{K}_m$, the sequence $(\overline{\Sigma}_m^+)$ also converges to the future side of $H$. However, by compactness, we may suppose that $N_m$ converges towards some limit $N_\infty$, say. In particular, this vector is a supporting normal at the origin to the limit of $(\overline{\Sigma}_m^+)$, that is, the future side of $H$. This is absurd, since $H$ is null, and the result follows.\qed
\proclaim{Lemma \nextprocno}
\noindent There exists a spacelike, locally strictly convex hypersurface $\Sigma_\infty$ towards which $\Sigma_m$ subconverges in the $C^\infty_\oploc$ sense.
\endproclaim
\proclabel{CompactnessForHypersurfaces}
\proof Upon rescaling, we may suppose that $k_m=-1/3$ for all $m$. Furthermore, by Lemma \procref{UniformlySpacelike}, the sequence $(\Sigma_m)$ is uniformly spacelike over every compact set. Since the scalar curvature is an elliptic curvature function (c.f. \cite{CaffNirSprV} and \cite{SmiRosDT}), it now suffices to show that for every compact subset $X$ of $\Bbb{R}^{d,1}$, there exists $B>0$ such that, for all $m$, the shape operator $A_m$ of $\Sigma_m$ satisfies
$$
\frac{1}{B} \leq A_m(x) \leq B,
$$
for all $x\in\Sigma_m\minter X$. Furthermore since, for all $m$, $\Sigma_m$ has constant scalar curvature equal to $(-1/3)$, it suffices to prove the lower bound. Now suppose the contrary, so that there exists a sequence of points $(x_m)$ contained within some compact set $X$, say, and a sequence $(\lambda_m)$ of positive real numbers converging to zero such that, for all $m$, $x_m$ is an element of $\Sigma_m$ and $\lambda_m$ is an eigenvalue of the shape operator of $\Sigma_m$ at this point.
\par
For all $m$, let $N_m$ be the future-oriented, unit, normal vector field over $\Sigma_m$, and define
$$\eqalign{
\hat{\Sigma}_m &:= \left\{ (x,N_m(x))\ |\ x\in\Sigma_m\right\},\ \text{and}\cr
\hat{x}_m &:= (x_m,N_m(x)).\cr}
$$
By Lemma \procref{Dictionary}, for all $m$, $\hat{\Sigma}_m$ is a complete, positive, special legendrian submanifold of $M:=\opU^+\Bbb{R}^{3,1}$ with refined special lagrangian angle equal to $\pi/2$. Furthermore, by Lemma \procref{UniformlySpacelike}, the sequence $(\hat{x}_m)$ is contained within a compact subset of $M$. It follows by Theorem \procref{CompactnessPSLC} that there exists a complete, pointed, positive, special legendrian submanifold $(\hat{\Sigma}_\infty,\hat{x}_\infty)$ towards which $(\hat{\Sigma}_m,\hat{x}_m)$ subconverges.
\par
Since the least eigenvalue of the shape operator of $\Sigma_m$ at the point $x_m$ tends to zero, it follows by Theorem \procref{TotalDegeneracy} that $\hat{\Sigma}_\infty$ is a curtain submanifold. In particular, by Theorem \procref{CurtainSubmanifolds} its projection onto $\Bbb{R}^{3,1}$ is foliated by complete, spacelike geodesics, which must all be contained in $K_\infty$. This is absurd, since the intersection of $K_\infty$ with any spacelike hyperplane is compact. The result follows.\qed
\newsubhead{Existence and uniqueness}[ExistenceAndUniquenss]
We now prove the main result of this paper.
\proclaim{Theorem \nextprocno}
\noindent Let $\Gamma\subseteq\opPSO(3,1)$ be a kleinian subgroup. Let $\tau:\Gamma\rightarrow\Bbb{R}^{d,1}$ be an $\Bbb{R}^{3,1}$-cocycle. Let $K:=K(\tau)$ be Bonsante's domain for $\tau$. For all $k>0$, there exists a unique, spacelike, locally strictly convex hypersurface $\Sigma_k$ in $K$ which is of constant scalar curvature equal to $(-k^2)$ and which is invariant under the action of $\Gamma\ltimes\tau$. Furthermore, the family $(\Sigma_k)_{k>0}$ constitutes a smooth foliation of the interior of $K$.
\endproclaim
\proclabel{MainTheoremWithFoliations}
\proof We first prove existence. Thus, let $I\subseteq[0,1]$ be the set of all $t\in[0,1]$ such that the existence part of the result holds for the cocycle $t\tau$. By Lemmas \procref{Openness} and \procref{CompactnessForHypersurfaces}, $I$ is both open and closed. The fuchsian case described at the end of Section \subheadref{QuasifuchsianRepresentations} shows that $0$ is an element of $I$ so that, by connectedness, $1$ is also an element of $I$, and existence follows.
\par
In fact, it follows from Lemmas \procref{JacobiOperator} and \procref{CompactnessForHypersurfaces} that the space of such hypersurfaces is discrete and compact and is therefore finite. Furthermore, one can then show via an elementary degree theoretic argument (c.f. \cite{SmiRosDT}) that the number of such hypersurfaces is independent of both $k$ and $\tau$. However, uniqueness can also be proven more directly as follows. Let $\Sigma_k$ be a spacelike, locally strictly convex hypersurface in $K$ of constant scalar curvature equal to $(-k^2)$ which is invariant under the action of $\Gamma\ltimes\tau$. For $t\in\Bbb{R}$, define $\Sigma_{k,t}:=\Sigma_k + t e_4$, where $e_4$ here denotes the fourth canonical basis element of $\Bbb{R}^{3,1}$. That is, $\Sigma_{k,t}$ is obtained by translating $\Sigma_k$ vertically upwards by a distance of $t$. Observe now that $\Sigma_k$ is asymptotic to $\partial K$ in the sense that the vertical distance between the two tends to $0$ at infinity. Now let $\Sigma_k'$ be another hypersurface with the same properties as $\Sigma_k$. Since $\Sigma_k'$ is also asymptotic to $\partial K$, there exists $t\in\Bbb{R}$ such that $\Sigma_k'$ is an interior tangent to $\Sigma_{k,t}$ at some point, and it follows by the strong geometric maximum principle that these two hypersurfaces coincide. Finally, since both $\Sigma_k$ and $\Sigma_k'$ are asymptotic to one another, it follows that $t=0$, and this proves uniqueness.
\par
We now prove that the family $(\Sigma_k)_{k>0}$ smoothly foliates the interior of $K$. Consider $k\in\Bbb{R}$. For $l$ sufficiently close to $k$, $\Sigma_l$ is the normal graph of some function $f_l$ over $\Sigma_k$. Consider now the partial derivative,
$$
g:=\frac{\partial}{\partial l}f_l|_{l=k}.
$$
By Lemma \procref{JacobiOperator},
$$
B^{ij}g_{;ij} = 1 + \phi g,
$$
for some positive-definite matrix $B^{ij}$ and some strictly positive function $\phi$. It follows by the maximum principle that $g$ is strictly negative and there therefore exists $\epsilon>0$ such that $(\Sigma_l)_{l\in]k-\epsilon,k+\epsilon[}$ smoothly foliates a neighbourhood of $\Sigma_k$. Since this holds for all $k$, it follows that $(\Sigma_k)_{k>0}$ smoothly foliates some open subset of $K$.
\par
It remains to show that this foliation covers the whole of the interior of $K$. However, by Lemma \procref{Position}, $\Sigma_k$ converges to $\partial K$ as $k$ tends to $+\infty$. On the other hand, for $t\in\Bbb{R}$ let $C_t$ be the future cone based on the point $te_4$, where $e_4$ again denotes the fourth canonical basis vector of $\Bbb{R}^{3,1}$. For $t$ sufficiently large and negative, $X$ is contained in $C_t$. Now let $(\Sigma_k')_{k>0}$ be the fuchsian foliation of $C_t$ constructed at the end of Section \subheadref{QuasifuchsianRepresentations}. By the geometric maximum principle, for all $k$, $\Sigma_k$ lies above $\Sigma_k'$. From this it follows that $(\Sigma_k)_{k>0}$ foliates the whole of the interior of $K$, and this completes the proof.\qed
\newhead{Bibliography}[Bibliography]
\medskip
{\leftskip = 5ex \parindent = -5ex
\leavevmode\hbox to 4ex{\hfil \cite{AnderssonBarbotBeguinZeghib}}\hskip 1ex{Andersson L., Barbot T., B\'eguin F., Zeghib A., Cosmological time versus CMC time in spacetimes of constant curvature, {\sl Asian Journal of Mathematics}, {\bf 16}, (2012), no. 1, 37--88}
\medskip
\leavevmode\hbox to 4ex{\hfil \cite{Apanasov}}\hskip 1ex{Apanasov B. N., Deformations of conformal structures on hyperbolic manifolds, {\sl J. Differential Geom.}, {\bf 35}, (1992), no. 1, 1--20}
\medskip
\leavevmode\hbox to 4ex{\hfil \cite{Barbot}}\hskip 1ex{Barbot T., Flat globally hyperbolic spacetimes,{\sl J. Geom. Phys.}, {\bf 53}, (2005), no.2, 123--165}
\medskip
\leavevmode\hbox to 4ex{\hfil \cite{BarbotBonsanteSchlenker}}\hskip 1ex{Barbot T., Bonsante F., Schlenker J.-M., Collisions of particles in locally AdS spacetimes I. Local description and global examples, {\sl Comm. Math. Phys.}, {\bf 308}, (2011), no. 1, 147--200}
\medskip
\leavevmode\hbox to 4ex{\hfil \cite{BarbotBeguinZeghibI}}\hskip 1ex{Barbot T., B\'eguin F., Zeghib A., Prescribing Gauss curvature of surfaces in 3-dimen\-sional spacetimes, Application to the Minkowski problem in the Minkowski space, {\sl Ann. Instit. Fourier.}, {\bf 61}, (2011), no. 2, 511--591}
\medskip
\leavevmode\hbox to 4ex{\hfil \cite{SanchezI}}\hskip 1ex{Bernal A. N., S\'anchez M., On smooth Cauchy hypersurfaces and Geroch’s splitting theorem, {\sl Comm. Math. Phys.}, {\bf 243}, (2003), no. 3, 461--470}
\medskip
\leavevmode\hbox to 4ex{\hfil \cite{SanchezII}}\hskip 1ex{Bernal A. N., S\'anchez M., Globally hyperbolic spacetimes can be defined as ‘causal’ instead of ‘strongly causal’, {\sl Classical Quantum Gravity}, {\bf 24}, (2007), no. 3, 745--749}
\medskip
\leavevmode\hbox to 4ex{\hfil \cite{Bonsante}}\hskip 1ex{Bonsante F., Flat spacetimes with compact hyperbolic Cauchy surfaces, {\sl J. Differential Geom.}, {\bf 69}, (2005), no. 3, 441--521}
\medskip
\leavevmode\hbox to 4ex{\hfil \cite{BonsanteFillastre}}\hskip 1 ex{Bonsante F., Fillastre F., The equivariant Minkowski problem in Minkowski space, to appear in {\sl Ann. Inst. Fourier}}
\medskip
\leavevmode\hbox to 4ex{\hfil \cite{BonsanteMondelloSchlenkerI}}\hskip 1ex{Bonsante F., Mondello G., Schlenker J.-M., A cyclic extension of the earthquake flow, {\sl Geometry \& Topology}, {\bf 17}, (2013), 157--234}
\medskip
\leavevmode\hbox to 4ex{\hfil \cite{BonsanteMondelloSchlenkerII}}\hskip 1ex{Bonsante F., Mondello G., Schlenker J.-M., A cyclic extension of the earthquake flow II, {\sl Annales scientifiques de l'ENS}, {\bf 48}, no. 4, (2015), 811--859}
\medskip
\leavevmode\hbox to 4ex{\hfil \cite{CaffNirSprV}}\hskip 1ex{Caffarelli L., Nirenberg L., Spruck J., Nonlinear second-order elliptic equations. V. The Dirichlet problem for Weingarten hypersurfaces, {\sl Comm. Pure Appl. Math.}, {\bf 41}, (1988), no. 1, 47--70}
\medskip
\leavevmode\hbox to 4ex{\hfil \cite{Carlip}}\hskip 1ex{Carlip S., {\sl Quantum gravity in $2+1$ dimensions}, Cambridge Monographs of Mathematical Physics, Cambridge University Press, Cambridge, (1998)}
\medskip
\leavevmode\hbox to 4ex{\hfil \cite{FillastreSmith}}\hskip 1ex{Fillastre F., Smith G., Group actions and scattering problems in Teichmueller theory, arXiv:1605.04563}
\medskip
\leavevmode\hbox to 4ex{\hfil \cite{GilbTrud}}\hskip 1ex{Gilbarg D., Trudinger N. S., {\sl Elliptic partical differential equations of second order}, Die Grundlehren der mathemathischen Wissenschaften, {\bf 224}, Springer-Verlag, Berlin, New York (1977)}
\medskip
\leavevmode\hbox to 4ex{\hfil \cite{Goldman}}\hskip 1ex{Goldman W. M., {\sl Discontinuous groups and the Euler class}, PhD Thesis, University of California, Berkeley, 1980, 138 pp.}
\medskip
\leavevmode\hbox to 4ex{\hfil \cite{JohnsonMillson}}\hskip 1ex{Johnson D., Millson J. J., Deformation spaces associated to compact hyperbolic manifolds, {\sl Discrete groups in geometry and analysis (New Haven, Conn., 1984)}, Progr. Math., vol. {\bf 67}, Birkh\"aauser Boston, Boston, MA, 1987, pp. 48--106}
\medskip
\leavevmode\hbox to 4ex{\hfil \cite{Kapovich}}\hskip 1ex{Kapovich M., Deformations of representations of discrete subgroups of $\opSO(3,1)$, {\sl Math. Ann.}, {\bf 299}, (1994), no. 2, 341--354}
\medskip
\leavevmode\hbox to 4ex{\hfil \cite{KapovichII}}\hskip 1ex{Kapovich M., Kleinian Groups in Higher Dimensions, in {\sl Geometry and Dynamics of Groups and Spaces}, Progress in Mathematics, {\bf 265}, 487--564}
\medskip
\leavevmode\hbox to 4ex{\hfil \cite{Kourouniotis}}\hskip 1ex{Kourouniotis C., Deformations of hyperbolic structures, {\sl Math. Proc. Cambridge Philos. Soc.}, {\bf 98}, (1985), no. 2, 247--261}
\medskip
\leavevmode\hbox to 4ex{\hfil \cite{KulkPink}}\hskip 1ex{Kulkarni R.S., Pinkall U., A canonical metric for M\"obius structures and its applications, {\sl Math. Z.}, {\bf 216}, (1994), no.1, 89--129}
\medskip
\leavevmode\hbox to 4ex{\hfil \cite{LabA}}\hskip 1ex{Labourie F., Un lemme de Morse pour les surfaces convexes, {\sl Invent. Math.}, {\bf 141}, (2000), no. 2, 239--297}
\medskip
\leavevmode\hbox to 4ex{\hfil \cite{Mess}}\hskip 1ex{Mess G., Lorentz spacetimes of constant curvature, {\sl Geom. Dedicata}, {\bf 126}, (2007), 3--45}
\medskip
\leavevmode\hbox to 4ex{\hfil \cite{Scannell}}\hskip 1ex{Scannell K. P., Infinitesimal deformations of some $\opSO(3,1)$ lattices, {\sl Pacific J. Math.}, {\bf 194}, (2000), no. 2, 455--464}
\medskip
\leavevmode\hbox to 4ex{\hfil \cite{SmiRosDT}}\hskip 1ex{Rosenberg H., Smith G., Degree Theory of Immersed Hypersurfaces, to appear in {\sl Mem. Amer. Math. Soc.}}
\medskip
\leavevmode\hbox to 4ex{\hfil \cite{SmiFCS}}\hskip 1ex{Smith G., Moduli of Flat Conformal Structures of Hyperbolic Type, {\sl Geom. Dedicata}, {\bf 154}, (2011), no. 1, 47--80}
\medskip
\leavevmode\hbox to 4ex{\hfil \cite{SmiSLC}}\hskip 1ex{Smith G., Special Lagrangian curvature, Math. Annalen, {\bf 335}, (2013), no. 1, 57--95}
\medskip
\leavevmode\hbox to 4ex{\hfil \cite{SpruckXiao}}\hskip 1ex{Spruck J., Xiao L., Convex Spacelike Hypersurfaces of Constant Curvature in de Sitter Space, {\sl Discrete Contin. Dyn. Syst. Ser. B}, {\bf 17}, no.6, (2012), 2225-–2242}
\par}
%
%
%
%
\enddocument